\documentclass[10pt,a4paper]{amsart}

\usepackage{amsfonts,amssymb,latexsym}

\def\LaTeX{L\kern-.36em\raise.3ex\hbox{a}\kern-.15em
    T\kern-.1667em\lower.7ex\hbox{E}\kern-.125emX}

\newcommand{\ka}[6]{{\ent{\begin{array}{ccccccccccc}
\ &\  &\     &\     &\     &  {\ensuremath{#1}^2}  &\     &\     &\
&\  &\ \\
\ &\  &\     &\     &  {\ensuremath{#2}}.{\ensuremath{#1}} &  0  &
{\ensuremath{#2}}.{\ensuremath{#1}} &\     &\     &\  &\ \\
\ &\  &\     &  {\ensuremath{#3}}.{\ensuremath{#1}} &   0
&{\ensuremath{#2}^2}&   0  &  {\ensuremath{#3}}.{\ensuremath{#1}} &\
&\  &\ \\
\ &\  &  {\ensuremath{#4}}.{\ensuremath{#1}} &   0
&{\ensuremath{#2}}.{\ensuremath{#3}}&  0
&{\ensuremath{#2}}.{\ensuremath{#3}}&   0  &
{\ensuremath{#4}}.{\ensuremath{#1}} &\  &\ \\
\ &{\ensuremath{#5}}.{\ensuremath{#1}}&   0
&{\ensuremath{#2}}.{\ensuremath{#4}}&   0  &{\ensuremath{#3}^2}&   0
&{\ensuremath{#2}}.{\ensuremath{#4}}&   0
&{\ensuremath{#5}}.{\ensuremath{#1}}&\ \\
{\ensuremath{#6}}.{\ensuremath{#1}} & 0
&{\ensuremath{#2}}.{\ensuremath{#5}}&   0
&{\ensuremath{#3}}.{\ensuremath{#4}}&  0
&{\ensuremath{#3}}.{\ensuremath{#4}}&   0
&{\ensuremath{#2}}.{\ensuremath{#5}}& 0
&{\ensuremath{#6}}.{\ensuremath{#1}} \\
\ &{\ensuremath{#2}}.{\ensuremath{#6}}&   0
&{\ensuremath{#3}}.{\ensuremath{#5}}&   0  &{\ensuremath{#4}^2}&   0
&{\ensuremath{#3}}.{\ensuremath{#5}}&   0
&{\ensuremath{#2}}.{\ensuremath{#6}}&\ \\
\ &\  &  {\ensuremath{#3}}.{\ensuremath{#6}} &   0
&{\ensuremath{#4}}.{\ensuremath{#5}}&  0
&{\ensuremath{#4}}.{\ensuremath{#5}}&  0  &
{\ensuremath{#3}}.{\ensuremath{#6}} &\  &\ \\
\ &\  &\     &  {\ensuremath{#4}}.{\ensuremath{#6}} &   0
&{\ensuremath{#5}^2}&  0  &  {\ensuremath{#4}}.{\ensuremath{#6}} &\
&\  &\ \\
\ &\  &\     &\     &  {\ensuremath{#5}}.{\ensuremath{#6}} &  0  &
{\ensuremath{#5}}.{\ensuremath{#6}} &\     &\     &\  &\ \\
\ &\  &\     &\     &\     &  {\ensuremath{#6}^2}  &\     &\     &\
&\  &\ \\
\end{array}}
}}

\def\aa{{\mathcal A}}
\def\bb{{\mathcal B}}

\def\ee{{\mathcal E}}
\def\ff{{\mathcal F}}

\def\lll{{\mathcal L}}

\def\pp{{\mathcal P}}
\def\qq{{\mathcal Q}}

\def\ss{{\mathcal S}}

\def\sN{{\ss(N)}}

\def\tq{\ :\ }

\def\eps{\varepsilon}
\def\dst{\displaystyle}

\def\supp{{\mathrm{supp}\,}}

\newcommand{\Lt}{L^2(\R)}

\newcommand{\LT}{L^2(\T)}

\newcommand{\Ll}{L_\Lambda^2(\T)}

\newcommand{\Ol}{\overline}

\newcommand{\Be}{\begin{equation}}
\newcommand{\Ee}{\end{equation}}

\def\udots{\raisebox{0.35mm}{\mbox{.}\raisebox{1.1mm}{\mbox{\hskip0.2mm.}}\raise
box{2.2mm}{\mbox{\hskip0.3mm.}}}}

\def\C{{\mathbb{C}}}

\def\hil{{\mathbb{H}}}

\def\N{{\mathbb{N}}}

\def\R{{\mathbb{R}}}

\def\T{{\mathbb{T}}}

\def\Z{{\mathbb{Z}}}
\def\la{\lambda}
\def\La{\Lambda}
\def\Lf{\Lambda_f}
\def\Lg{\Lambda_g}
\def\al{\alpha}
\def\dt{\delta}
\def\hP{(\hat{\rm P})}
\def\bc{{\bf c}}
\def\hf{\hat{f}}
\def\hg{\hat{g}}

\def\tK{\widetilde{K}}
\def\th{\hat{h}}
\newcommand{\Pl}{\rm{(}\widehat{P}_{\Lambda}\rm{)}}

\newcommand{\norm}[1]{{\left\|{#1}\right\|}}
\newcommand{\ent}[1]{{\left[{#1}\right]}}
\newcommand{\abs}[1]{{\left|{#1}\right|}}
\newcommand{\scal}[1]{{\left\langle{#1}\right\rangle}}
\newcommand{\set}[1]{{\left\{{#1}\right\}}}
\newenvironment{notation}[1][]{\vskip1pt\noindent\rm\textbf{Notation}.\
 }{\rm \vskip 1pt}
\newenvironment{remark}[1][]{\vskip1pt\noindent\rm\textit{Remark}:\
}{\rm\vskip1pt}
\newenvironment{defi}[1][]{\vskip3pt\noindent\it\textbf{Definition.}\
}{\rm\vskip3pt}
\newenvironment{example}[1][]{\vskip1pt\noindent\rm\textit{Example\
#1}:\ }{\rm\vskip1pt}

\newtheorem{problem}{Problem}

\newtheorem{lem}{Lemma}[section]
\newtheorem{prop}[lem]{Proposition}
\newtheorem{theo}[lem]{Theorem}
\newtheorem{cor}[lem]{Corollary}

\newtheorem{remit}[lem]{Remark}

\newenvironment{lemma}[1][]{\lem\sl}{\rm\vskip3pt}
\newenvironment{proposition}[1][]{\prop\sl}{\rm\vskip3pt}
\newenvironment{theorem}[1][]{\theo\sl}{\rm\vskip3pt}
\newenvironment{corollary}[1][]{\cor\sl}{\rm\vskip3pt}

\newenvironment{prob}[1][]{\problem\sl}{\rm\vskip3pt}
\newenvironment{remarknum}[1][]{\remit\sl}{\rm\vskip3pt}

\def\drap{\simeq}
\def\tdrap{\equiv}
\def\ntdrap{\sim}

\begin{document}
\title{Discrete radar ambiguity problems}
\author{Aline BONAMI, Gustavo GARRIG\'OS \& Philippe JAMING}
\address{A.~B. \& P.~J.\,: MAPMO\\ Universit\'e d'Orl\'eans\\ BP 6759\\ F 45067 ORLEANS Cedex 2\\
FRANCE}
\email{Aline.Bonami@univ-orleans.fr,\ Philippe.Jaming@univ-orleans.fr}
\address{G.~G.~: Universidad Aut\'onoma de Madrid \\
Departamento de Matem\'aticas C-XV\\
Ciudad Universitaria Cantoblanco\\
28049 Madrid\\
SPAIN}
\email{gustavo.garrigos@uam.es}


\date{}

\thanks{Research partially financed by : {\it European Commission}
Harmonic Analysis and Related Problems 2002-2006 IHP Network
(Contract Number: HPRN-CT-2001-00273 - HARP).}

\begin{abstract}
In this paper, we pursue the study of the radar ambiguity problem
started in \cite{Ja,GJP}. More precisely, for a given function $u$
we ask for all functions $v$ (called \emph{ambiguity partners}) such that the ambiguity functions of $u$
and $v$ have same modulus. In some cases, $v$ may be given by some elementary transformation
of $u$ and is then called a \emph{trivial partner} of $u$ otherwise we call it a \emph{strange partner}.
Our focus here is on two discrete versions of the problem.

For the first one, we restrict the problem to functions $u$ of the Hermite class, $u=P(x)e^{-x^2/2}$,
thus reducing it to an algebraic problem on polynomials. Up to some mild restriction
satisfied by quasi-all and almost-all polynomials, we show that such a function has only trivial partners.

The second discretization, restricting the problem to pulse type signals,
reduces to a combinatorial problem on matrices of a special form.
We then exploit this to obtain new examples of functions that
have only trivial partners. In particular, we show that
most pulse type signals have only trivial partners.

Finally, we clarify the notion of \emph{trivial partner},
showing that most previous counterexamples are still trivial in some restricted sense.
\end{abstract}

\maketitle

\tableofcontents
\newpage
\section{Introduction}
\label{intro}

Phase retrieval problems arise naturally in the applied study of
signals \cite{Wa,Hu,KST,tes}...
They are based on the ambiguity for the phase choice
in a signal with fixed frequency amplitude.
To be more precise, let us denote the Fourier transform
of $u\in L^1(\R)$ (with the usual extension to $L^2(\R)$) by $\ff$:
$$
\ff u(\xi)=\int_\R u(x)\,e^{i x\xi}\,dx,\quad\xi\in\R.
$$
The {\it phase retrieval problem} then amounts to solving
the following:

\begin{prob}[Phase retrieval] Given $u\in L^2(\R)$, find all $v\in L^2(\R)$ such
 that for all
$x\in\R$,
\begin{equation}
\abs{\ff u(x)}=\abs{\ff v(x)}.
\label{prp}
\end{equation}
\end{prob}
This problem admits always the trivial solutions $v(x)=c\,u(x-\alpha)$ and
 $v(x)=c\,\overline{u(-x-\alpha)}$,
where $|c|=1$ and $\alpha\in\R$.

In applied problems, one may usually further restrict the class of functions to
 which $u$ and $v$ should
belong. A typical example would be to ask for $u$ and $v$ to be compactly
 supported.
In this case, there are usually many non trivial solutions and a complete
 description of them
is available in terms of the zeros of the holomorphic functions
$\ff u$ and $\ff v$ (see \cite{Wa,Ro,Ja} for a complete
description of these
 solutions). For further information on phase retrieval problems,
we refer to these articles as well as \cite{JK}, the surveys \cite{KST,Mi},
the book \cite{Hu} and references therein.

In this paper we shall deal with a different, although closely
related type of phase retrieval problem, having its origin in the
analysis of radar signals. Following Woodward \cite{Wo}, a radar antenna emits a signal
$u\in\Lt$ that is reflected by a target and modified by Doppler
effect. It then returns to the antenna where it is correlated by
the emitted signal, so that, under certain physical conditions,
the radar measures the quantity:
\begin{equation}
A(u)(x,y) =  \dst\int_\R u\left(t\right)
\overline{u\left(t-x\right)}e^{i yt}dt,
\quad x,y\in\R,
\label{raf}
\end{equation}
and $A(u)$ is called the {\sl radar ambiguity function of $u$}. As
usually happens, receivers are not able to read the phase, but
only the amplitude $|A(u)(x,y)|$, giving rise to the following
{\it radar ambiguity problem}:

\begin{prob}[Radar ambiguity]\label{prob3}
Let $u\in \Lt$, then find all $v\in\Lt$ such that
\begin{equation}
|A(u)(x,y)|=|A(v)(x,y)|,\quad x,y\in\R.
\label{rap}
\end{equation}
\end{prob}
Note that, for each $x\in\R$,
$A(u)(x,\cdot)=\ff\ent{u\left(\cdot\right)
{\overline{u}}\left(\cdot-x\right)}$, so
that equation (\ref{rap}) is actually a family of phase retrieval
problems as described in (\ref{prp}). Two functions $u$ and $v$
satisfying (\ref{rap}) are said to be
{\it (radar) ambiguity partners}.
The reader may find a comprehensive historical introduction and further
references to this problem in \cite{Ja}. Properties of $A(u)$ that we
may use in this paper can all be found there and in \cite{AT,Wi}
(note that we slightly change the normalization for $A(u)$ from \cite{AT}).

It is not difficult to verify that trivial solutions to
the equation in (\ref{rap}) are given by:
\begin{equation}
v(t)=c\,e^{i\beta t}\,u(t-\alpha)\quad\mbox{and}\quad
v(t)=c\,e^{-i\beta t}\,u(-t-\alpha),\quad|c|=1,\,\,\alpha,\beta\in\R.
\label{hei}
\end{equation}
The first set of solutions corresponds to a unitary
representation of the Heisenberg group,
while the second is just a composition with the isometry
$Z f(t)=f(-t)$.
So, following \cite{Ja}, we say that $u$ and $v$ are
{\it trivial partners}
 when they satisfy (\ref{hei}).
If $u$ and $v$ are ambiguity partners that are not trivial
partners, we will say that they are {\it strange partners} and in
\cite{dB,GJP,Ja}, examples of signals having strange partners are
given. In the opposite direction, there exist signals for which
every ambiguity partner is trivial.

The aim of this paper is to get some insight on which functions
may or may not have strange partners. To tackle this problem we
appeal to two different discrete (finite dimensional) versions of
the problem, both being also of practical interest.

The first discretization is the restriction of the problem to
Hermite functions, that is to functions of the form
$P(t)e^{-t^2/2}$ where $P$ is a polynomial. There are several
reasons for this: first it was proposed by Wilcox in his pioneering
paper \cite{Wi}, since it is a dense class of functions which
are best localized in the time-frequency plane and are thus well
adapted for numerical analysis. Second, in some sense this class is
``extremal'' for the \emph{uncertainty principle}, so
one can show that all solutions to Problem \ref{prob3}
are necessarily Hermite functions $v(t)=Q(t)e^{-t^2/2}$ for some polynomial $Q$ (except perhaps for trivial transformations;
{\it see} \cite{dB} or Lemma \ref{uncert} below).
Finally, Hermite functions are of
theoretical importance for the problem considered. Indeed,
Bueckner \cite{Bu} associated to each function $u\in L^2(\R)$ an
Hilbert-Schmidt operator $K_u$ in a way that finding all solutions
for the ambiguity problem for $u$ amounts to finding all functions
$v$ such that $K_u^*K_u=K_v^*K_v$. He then proved that $K_u$ is of
finite rank if and only if $u$ is a Hermite function.
Moreover, the following conjecture was proposed:

\medskip

\noindent{\bf Conjecture \cite{Bu}.} {\sl If $u$ is a Hermite function, then
$u$ has only trivial partners.}

\medskip

Indeed, Bueckner was considering the {\it bilinear} version of (\ref{rap}):
\begin{equation}
\abs{A(u_1,u_2)}=\abs{A(v_1,v_2)}
\label{bilradprob}
\end{equation}
where $A(u_1,u_2)$ is the bilinear functional associated with $A(u)$.
He proved that for almost every couple of functions of the form
$$
(u_1,u_2)=(P_1(x)e^{-x^2/2},P_2(x)e^{-x^2/2})
$$
($P_1,P_2$ polynomials), the solutions to (\ref{bilradprob}) are trivial
partners of $(u_1,u_2)$. However, his techniques depend on a certain criterion
that excludes the quadratic case, and hence do not say anything about Problem
\ref{prob3}.

In this paper we will prove, using a simple algebraic approach, the following result about ambiguity
partners of Hermite functions:
\vskip6pt
\noindent{\bf Theorem A.} {\sl For almost all and quasi-all polynomials $P$,
the function $u(x)=P(x)e^{-x^2/2}$
has only trivial partners.}
\vskip6pt
Here almost-all (respectively quasi-all)
refers to Lebesgue measure (respectively Baire category) when one identifies the set of polynomials of fixed degree $n$
with $\C^{n+1}$.

The problem has also been considered by deBuda \cite{dB}, who obtained some partial results
in an unpublished report which unfortunately are not always complete. Although our approach
shares some common features with his, it is essentially distinct
as we introduce a new argument by using the fact that $A(u)$ has
some factorization if $u$ has non-trivial partners. Some technical
difficulties remain as our use of Bezout's theorem forces us to
assume that some polynomial associated to $u$ has only simple
non-symmetric zeros in order to prove that $u$ has only trivial
partners.

\medskip

The second class of functions we consider is the restriction to
signals of {\it pulse type}:
\begin{equation}
u(t)=\sum_{j=-\infty}^\infty a_jH(t-j),\quad x\in\R,
\label{puls}
\end{equation}
where $H\in\Lt$ has $\supp H\subset\,\ent{0,\frac{1}{2}}$, and $\{a_j\}_{j\in\Z}$
is a finite sequence of complex numbers. This class of functions is very common in radar signal design
({\it see e.g.} \cite[p 285]{van}). It also leads naturally to a discretization
of Problem \ref{prob3}. Indeed, a simple computation shows that, for all $k\in\Z$, $y\in\R$ and
$k-\frac12\leq x\leq k+\frac12$, one has:
\begin{equation}
A(u)(x,y)=\,\biggl(\sum_{j\in\Z}a_j{\Ol{ a_{j-k}}}e^{i jy}
\biggr)
\,A(H)(x-k,y).
\label{pulse}
\end{equation}
This following \emph{discrete ambiguity problem} was proposed in \cite{GJP}:

\begin{problem}[Discrete Radar Ambiguity Problem]
Given $a=\{a_j\}\in\ell^2(\Z)$, find all sequences $\,
b\in\ell^2(\Z)$ such that, for every $k\in\Z$ and  $y\in\R$,
\begin{equation}
\abs{\aa(a)(k,y)}=\abs{\aa(b)(k,y)},
\label{pbdisc}
\end{equation}
where
$$
\aa(a)(k,y)=\sum_{j\in\Z}a_j{\Ol{ a_{j-k}}}e^{i jy}.
$$
\label{problemP}
\end{problem}

Again, a sequence $b$, solution to (\ref{pbdisc}), is called an
\emph{ambiguity partner of $a$}. It is easy to see that trivial solutions
to (\ref{pbdisc})are given by
$$
b_j=ce^{i\beta j}a_{j-k}\quad\mbox{and}\quad
b_j=ce^{i\beta j}a_{-j-k},\quad|c|=1,\ \beta\in\R,\ k\in\Z.
$$
Such solutions are again called \emph{trivial partners of $a$} and
solutions that are not of this type are called \emph{strange
partners}. The main result of \cite{GJP} shows that a finite
sequence $a=\{a_j\}\in\C^{d+1}$ has only trivial partners, except
perhaps for $a$'s in a semialgebraic set of real codimension $1$
in $\C^{d+1}$ ({\it see} Theorem \ref{th:pulserare} below). This
was done by adapting Bueckner's method to the Discrete Radar
Ambiguity Problem, and then adapting a careful analysis to the
obtained combinatorial equation of matrices. The form of these
matrices was also exploited to produce new constructions of
non-trivial solutions in the exceptional set. A few other points
about such constructions, which were only announced in \cite{GJP},
are proven here in full detail ({\it see} Section \ref{sec:4.3}).

It was not investigated, however, how to translate these discrete
results into uniqueness statements for the general ambiguity
problem {\it i.e.} to Problem \ref{prob3}. This step is now
different from the corresponding one for Hermite functions, since
the class of pulse type signals is not extremal for the
uncertainty principle. In this paper, we introduce new techniques
for this class based on complex analysis and distribution theory,
which allows us to prove the following theorem:

\medskip

{\bf Theorem B.} {\sl Let $0<\eta\leq\frac{1}{3}$, and let
$a=(a_0,a_1,\ldots,a_N)\in\C^{N+1}$ that has only trivial
partners. Then the pulse type signal
$$
u(t)=\sum_{j=0}^Na_j\chi_{[j,j+\eta]}(t).
$$
has only trivial partners.}

\medskip

We do not know whether the condition $\eta\leq1/3$ is optimal. It was essential
in the proof to ensure that $v$ is also of pulse type.

\medskip

Next, we clarify the notion of \emph{trivial solutions}. There are
numerous phase retrieval problems in the literature and we think
that a natural definition of a trivial solution is to be a linear
or anti-linear operator that associates to each function a
solution of the given phase retrieval problem. Using Theorem A, we
will show that those trivial solutions described in Equation
(\ref{hei}) are indeed the only trivial solutions in the previous
sense:

\medskip

{\bf Theorem C.} {\sl The only linear (or anti-linear) bounded transformations $
T:\ L^2(\R)\to L^2(\R)$ so that
$$
|A(Tu)(x,y)|=|A(u)(x,y)|,\quad \mbox {for all }u\in L^2(\R)
$$
are those described in (\ref{hei}).}

\medskip

We do not know of an earlier proof of that simple fact. This
theorem is also reminiscent of Wigner's Unitary-Antiunitary
Theorem ({\see e.g.} \cite{LM,Ra,Mol}) which can be stated
as follows. Let $T$ be an operator $T$ on a Hilbert space $H$ and
assume that $T$ preserves the modulus of the scalar product:
$$
|\scal{Tx,Ty}|=|\scal{x,y}|\quad \mbox{for all}\ x,y\in H.
$$
Then $T$ is of the form $Tx=\omega(x)Ux$ where $\omega$
is a scalar valued function on $H$ such that $|\omega(x)|=1$
and $U$ is either unitary or anti-unitary operator on $H$.
Here we are in a slightly different situation and
Wigner's theorem can not be applied. It does nevertheless ask
whether the (anti)linearity assumption in Theorem C may be removed.
\medskip

Finally, we also consider a further restriction of the Discrete
Ambiguity Problem by considering sequences in $\ell^2(\Lambda)$
for some $\lambda\subset\Z$.
This is natural since most of the known examples of signals with strange
partners are of the form
$$
u(t)=\sum_{j\in\Lambda} c_j \chi_{[0,T]+j},
$$
at least when $\Lambda$ has ``enough gaps'' ({\it see e.g.} \cite{Ja}).
Indeed, partners of $u(t)$ can be easily obtained by multiplying
each $c_j$ by a unimodular constant $\exp(i\omega_j)$.
Here we clarify the nature of these ``gaps'' in terms of arithmetic conditions
which appear in the classical theory of trigonometric series with gaps.
More precisely, we assume that $\Lambda$ is a $B_2$ or a $B_3$-set
({\it see} Remark \ref{two} for precise definitions). In particular, and as a
consequence of our results we obtain the following.

\medskip

\noindent{\bf Theorem D.} {\sl Let $\dst u(t)=\sum_{j\in \Lambda} c_j \chi_{[0,T]+j}$.
Then, if $\Lambda$ is a $B_2$-set, then for all real $\omega_j$,
$$
v(t)=\sum_{j\in \Lambda} e^{i\omega_j} c_j \chi_{[0,T]+j}
$$
is a partner of $u(t)$. Moreover, if $\Lambda$ is a finite $B_3$-set
these are all partners of $u$.}

\medskip

Nevertheless, recall from \cite{GJP} or Forumla (\ref{ex:triv}) below,
that already when $\Lambda=\{0,1,2,3\}$ there exist
exceptional cases when strange solutions cannot be classified in
terms of gaps.

\medskip

The article is organized as follows. In the next section, we
concentrate on the continuous problem for Hermite functions, and
we prove Theorem A. The following section is devoted to the
characterization of trivial solutions, both in the discrete case
and in the continuous case. The last section is devoted to the
case of pulse type signals.
We start by proving Theorem B and conclude by
recalling and completing the main results of \cite{GJP}.

\section{The ambiguity problem for Hermite functions}

We now prove Theorem A. We will need a certain number of steps in
the proof. The two first ones are mainly due to DeBuda, \cite{dB}
and \cite{dB2}. In particular, De Buda has established the
stability of the class of Hermite signals for the ambiguity
problem using an elementary proof (which is not complete in
\cite{dB}). It can also be obtained as a consequence of the
\emph{uncertainty principle for ambiguity functions}, as it is
mentioned in \cite{BDJ}.

\subsection{Stability of Hermite functions for the Ambiguity problem}

\begin{lemma}
\label{uncert} Let $u(t)=P(t)e^{-{t^2}/2}$, where $P(t)$ is a
polynomial. Then, except perhaps for a trivial transformation,
every ambiguity partner $v$ of $u$ is of the form
$v(t)=Q(t)e^{-{t^2}/2}$, where $Q(t)$ is a polynomial with
$\deg{P}=\deg{Q}$.
\end{lemma}

\begin{proof}
Using the fact
$\mathcal{F}(e^{-{t^2}/2})(\xi)=\sqrt{2\pi}\,e^{-{\xi^2}/2}$, an
elementary computation shows
$$
Au(x,y)=e^{-i\frac{xy}2}\,\widetilde{P}(x,y)\,e^{-\frac{x^2+y^2}4},
$$
where $\widetilde{P}(x,y)$ is a polynomial of 2 variables of total
degree $2\deg{P}$ (see, e.g., \cite[Theorem 7.2]{Wi} or
(\ref{eq:herlag}) below). Then,
$$
|Av(x,y)|^2=|Au(x,y)|^2=|\widetilde{P}(x,y)|^2\,e^{-\frac{x^2+y^2}2},
$$
so we can use the uncertainty principle in \cite[Prop. 6.2]{BDJ}
to conclude $v(t)=Q(t)e^{i\omega t}e^{-\frac{(t-a)^2}2}$, for a
polynomial $Q$ and two real constants $\omega,a$. We only need to
show that $\deg{Q}=\deg{P}$, but this follows easily from
$$
|Av(x,y)|=|\widetilde{Q}(x,y)|\,e^{-\frac{x^2+y^2}4}
=|\widetilde{P}(x,y)|\,e^{-\frac{x^2+y^2}4},
$$
and the fact
$2\deg{Q}=\deg{\widetilde{Q}}=\deg{\widetilde{P}}=2\deg{P}$.
\end{proof}

\subsection{Reformulation of the ambiguity problem as an algebraic problem}\
Let us first give some notation that we will use in this section.
\begin{notation}
We say that a polynomial is {\sl monic } when the coefficient of
its term of higher degree is equal to $1$.

For a polynomial $\Pi\in \mathbb{C}[Z]$, we will write $\Pi^*$ the
polynomial given by $\Pi^*(z)=\overline{\Pi(\bar z)}$.

For a polynomial $\Pi$ of degree $n$, that is, $\Pi\in
\mathbb{C}_n[Z]$,  we write $\check\Pi(z)=(-1)^n\Pi(-z)$. Note
that $(\Pi')\check{\ }=(\check\Pi)'$. We will thus write
unambiguously $\check\Pi'$. Remark also that $\check \Pi$ is monic
when $\Pi$ is.

For $\Pi,\Psi$ two polynomials, we write
\begin{equation}
\{\Pi,\Psi\}_-=\Pi\check\Psi-\check\Pi\Psi,
  \hspace{2cm}
\{\Pi,\Psi\}_+=\Pi\check\Psi+\check\Pi\Psi.
\label{crochets}
\end{equation}

\end{notation}

We shall prove that the ambiguity problem for Hermite functions is
equivalent to an algebraic problem, which we state now. For
$\pp\in \mathbb{C}_n[Z]$, we define its ambiguity polynomial as
the polynomial in two variables given by
$$
A_{\pp}(z, w):= \,\sum_{m=0}^n\, \frac{1}{m!}\,\pp^{(m)}(
z)\,\pp^{*\,(m)}( w).
$$
Note that $A_{\pp}=A_{\qq}$ if and only if there exists some
unimodular constant $c$ such that $\pp=c\qq$.

The ambiguity problem for Hermite functions will then be reduced
to the following one:
\medskip

\noindent {\bf The algebraic ambiguity problem.} {\sl For a given
polynomial $\pp$ of degree $n$, find all polynomials $\qq$ for
which one has the following identity:}
\begin{equation}
\label{algebraic}
  A_{\pp}(z,w)A_{\pp}(-z,-w)=A_{\qq}(z,w)A_{\qq}(-z,-w).
\end{equation}
 Again, our
question is the following: does there exist other partners than
the trivial ones, given by $c\pp$ and $c\check\pp$, with $c$ a
unimodular constant?

We first prove the equivalence between the two problems.

Let us denote by
$$
H_k(x)=(-1)^ke^{x^2}\frac{d^k}{dx^k}(e^{-x^2}),\quad
k=0,1,2,\ldots
$$
the Hermite polynomials. We recall that, with the normalizing
constant $\gamma_k=(\sqrt{\pi}2^kk!)^{\frac{1}{2}}$, the system
$$
\psi_k(x)=\frac{1}{\gamma_k}\,H_k(x)\,e^{-{x^2}/2},\quad
k=0,1,2,\ldots
$$
is an orthonormal basis of $L^2(\R)$, called the \emph{Hermite
basis} of $L^2$.

Let $\bb$ be the linear map on $\mathbb{C}[Z]$ defined by
$\bb(H_k)=2^{k/2}Z^k$ ({\it i.e.} $\bb$ is the Bargmann
transform).

The equivalence between the two problems is given by the following
lemma, which is essentially contained in \cite{dB2}.

\begin{lemma}
Let $P$ and $Q$ be two polynomials. Then $Pe^{-t^2/2}$ and
$Qe^{-t^2/2}$ are ambiguity partners if and only if $\bb(P)$ and
$\bb(Q)$ are partners for the algebraic ambiguity problem.
\end{lemma}

\begin{proof}
Consider the expansion of $P$ and $Q$ in terms of the basis of
Hermite polynomials
\begin{equation}
P=\sum_{j=0}^n\alpha_jH_j\mbox{ and }Q=\sum_{j=0}^n\beta_jH_j
\label{h.1}
\end{equation}
with $\alpha_n\not=0$, $\beta_n\not=0$.

First of all, an explicit computation gives the well-known formula
(see, e.g., \cite[Theorem 7.2]{Wi}):
\begin{equation}
\label{eq:herlag}
A(H_je^{-t^2/2},H_ke^{-t^2/2})(x,y)=\lll_{jk}(x/\sqrt{2},y/\sqrt{2})\,
e^{-(x^2+y^2)/4}\,e^{i\frac{xy}{2}}
\end{equation}
where $\lll_{j,k}$ is the Laguerre polynomial defined by
$$
\lll_{j,k}(x, y)=
\gamma_j\gamma_k\,\sqrt{\frac{k!}{j!}}(x+iy)^{j-k}\sum_{\ell=0}^k
\begin{pmatrix}j\\k-\ell\\
\end{pmatrix}\frac{(-1)^\ell}{\ell!}\,(x^2+y^2)^\ell,\quad\mbox{if}\quad j\geq k
$$
(and $\lll_{j,k}(x,y)=\lll_{k,j}(-x,y)$ if $j<k$). We can write
this formula in a unified way as:
\[
\lll_{j,k}(x,y)=\sqrt{\pi2^{j+k}}\,j!\,k!\,\sum_{m=0}^{j\wedge k}
\frac{(x+iy)^{j-m}(-x+iy)^{k-m}}{(j-m)!(k-m)!m!}.
\]
Thus, defining the new variable $z=x+iy$ we have
\begin{eqnarray*}
A(H_je^{-t^2/2},H_ke^{-t^2/2})(x,y) & = &
\sqrt{\pi}\,j!\,k!\,\Bigl(\,\sum_{m=0}^{j\wedge
k}\,\frac{2^{m}}{m!}\,
\frac{z^{j-m}(-\overline{z})^{k-m}}{(j-m)!(k-m)!}\,\Bigr)\,
e^{-|z|^2/4}\,e^{i\frac{xy}{2}}\\
& = & \sqrt{\pi}\,j!\,k!\,\Bigl(\,\sum_{m=0}^{j\wedge k
}\,\frac{2^{m}}{m!}\,\frac{\partial^m}{\partial t^m}
\left.\left(\frac{t^{j}}{j!}\right)\right|_{t=z}\,\frac{\partial^m}{\partial
t^m}
\left.\left(\frac{t^{k}}{k!}\right)\right|_{t=-\overline{z}}\Bigr)\,
e^{-|z|^2/4}\,e^{i\frac{xy}{2}}.
\end{eqnarray*}
So, calling $\pp:=\bb(P)=\sum_{j=0}^n\al_j2^{j/2}Z^j$, and using the
bilinearity of the operator $A$ we have
\[
|A(u)(x,y)|=\sqrt{\pi}\,\left|\,\sum_{m=0}^n\,
\frac{1}{m!}\,\pp^{(m)}(\sqrt 2z)\,{\overline{\pp^{(m)}(-\sqrt 2 z
)}}\right|\, e^{-|z|^2/4}.
\]
Calling $\qq:=\bb(Q)=\sum_{j=0}^n\beta_j 2^{j/2}Z^j$, the fact that
$u$ and $v$ are partners is equivalent to the identity
\begin{equation}
\label{eq:polfund} \left|\,\sum_{m=0}^n\, \frac{1}{m!}\,\pp^{(m)}(
z)\,{\overline{\pp^{(m)}(- z)}}\right|^2 =\left|\,\sum_{m=0}^n\,
\frac{1}{m!}\,\qq^{(m)}(z)\,{\overline{\qq^{(m)}(-z)}}\right|^2
\end{equation}
for all complex numbers $z$. Since two holomorphic polynomials in two
complex variables $z,w$ coincide when they coincide for
$z=-\overline w$, this is equivalent to the identity
\begin{eqnarray}
 \left(\,\sum_{m=0}^n\, \frac{1}{m!}\,\pp^{(m)}(
z)\,\pp^{*\,(m)}(w)\right)\left(\,\sum_{m=0}^n\,
\frac{1}{m!}\,\pp^{(m)}(
-z)\,\pp^{*\,(m)}(-w)\right)\hspace{4cm}\nonumber\\
\hspace{4cm}=\left(\,\sum_{m=0}^n\, \frac{1}{m!}\,\qq^{(m)}(
z)\,\qq^{*\,(m)}(w)\right)\left(\,\sum_{m=0}^n\,
\frac{1}{m!}\,\qq^{(m)}(
-z)\,\qq^{*\,(m)}(-w)\right).\label{eq:polfund2}
\end{eqnarray}
We recognize the algebraic ambiguity problem, which finishes the
proof of the lemma.
\end{proof}

\begin{remarknum}
\label{rem:normalplo} Note that the highest order coefficient in
(\ref{eq:polfund2}) is $|\alpha_n|^4=|\beta_n|^4$, so that
$|\beta_n|=|\alpha_n|$. Replacing $\qq$ by its trivial partner
$\widetilde{\qq}=\frac{\alpha_n}{\beta_n}\qq$, we may thus assume
that $\beta_n=\alpha_n$. Then, using the homogeneity of Equation
(\ref{eq:polfund2}), there is no loss of generality to assume that
$\beta_n=\alpha_n=1$.
\end{remarknum}

\subsection{Solution of the algebraic ambiguity problem in the generic case}\

\begin{defi}
By a generic polynomial $\pp$ we mean a polynomial that has only
simple roots and has no common root with $\check \pp$, that is, $\pp$
has only simple non-symmetric roots.
\end{defi}
Of course, almost all and quasi-all polynomials are generic.

\medskip

We will now prove the following theorem which implies Theorem A.

\begin{theorem}
Assume that the polynomial $\pp$ is generic and let $\qq$ be a
partner of $\pp$. Then $\qq$ is a trivial partner, that is, there
exists a unimodular constant $c$ such that either $\qq=c\pp$ or
$\qq=c\check\pp$. \label{generic}\end{theorem}

The proof is divided into two steps. In the first one, we will
directly use equation (\ref{eq:polfund2}) to get substantial
information on $\qq$. The second step will consist in exploiting
the factorization that $\aa(\pp)$ would have if $\qq$ was not a
trivial partner.

\begin{proof}[First step]
As explained in Remark \ref{rem:normalplo}, we can assume that
$\pp$ and $\qq$ are monic polynomials and write
$$
\pp:=Z^n+p_1Z^{n-1}+\cdots
+p_{n-1}Z+p_n,\hspace{1.5cm}\qq:=Z^n+q_1Z^{n-1}+\cdots
+q_{n-1}Z+q_n.
$$
Equation (\ref{eq:polfund2}) can as well be written $$
A_{\pp}A_{\check\pp}=A_{\qq}A_{\check\qq}.$$ Looking at
$A_{\pp}\in \mathbb{C}[Z,W]$ as a polynomial in $W$ with
coefficients in $\mathbb{C}[Z]$, we  can write
$$
A_{\pp}\equiv \pp W^n +\left(\bar p_1 \pp+n \pp'\right)W^{n-1}+
\left(\bar p_2 \pp+(n-1)\bar p_1\pp'+
\frac{n(n-1)}{2}\pp''\right)W^{n-2}
$$
modulo terms of smaller degree. Looking at the coefficient of
$W^{2n}$ in (\ref{eq:polfund2}), we get
\begin{equation}
\label{P=Q} \pp\check \pp=\qq\check\qq,
\end{equation}
which in particular implies that
\begin{equation}\label{derivate}
\pp'\check \pp+\pp\check \pp' =\qq'\check\qq+\qq\check\qq'
\end{equation}
and
$$
\pp''\check \pp+\pp\check \pp'' + 2\pp'\check
\pp'=\qq''\check\qq+\qq\check\qq''+ 2\qq'\check\qq'.
$$
Then, looking at the coefficient of $W^{2n-2}$ in
(\ref{eq:polfund2})\footnote{The coefficient of $W^{2n-1}$ only
leads to  Equation (\ref{derivate}).}, an elementary computation
which uses the previous identities leads to
\begin{equation}
n\pp'\check\pp'+\bar p_1\bigl(\pp\check \pp'-\check\pp\pp'\bigr)
=n\qq'\check\qq'+\bar q_1\bigl(\qq\check \qq'-\check\qq\qq'\bigr).
\label{P'=Q'}
\end{equation}
The highest order term in this equation gives $|q_1|=|p_1|$.

From (\ref{P=Q}) we deduce that there exist two monic polynomials
$A$ and $B$ such that
\begin{equation}
\label{fact}
  \pp=AB \quad\mbox{and}\quad\qq= A\check B.
\end{equation}
 Let us further write
$$
A:=Z^k+a_1Z^{k-1}+\cdots+a_k\quad\mbox{and}\quad
B:=Z^l+b_1Z^{l-1}+\cdots+a_l.
$$
Then $p_1=a_1+b_1$, $q_1=a_1-b_1$ and $|q_1|=|p_1|$ is equivalent to
\begin{equation}
\label{orth}
  a_1\bar b_1+\bar a_1 b_1=0.
\end{equation}
 These relations, written for all possible decompositions of $\pp$
as a product $AB$, is sufficient to prove that the set of
coefficients $p_1,\cdots, p_n$ is contained in a real analytic
variety of codimension $1$ in $\C^n$, and imply Theorem A. We will
not give details for this reduction since we have more
information, as stated in Theorem \ref{generic}.

Note that, using the notations defined by (\ref{crochets}), (\ref{P'=Q'}) may as well be written as
\begin{equation}
 \label{concl}
2\bar a_1A\check A\{B',B\}_-+2\bar b_1B\check B\{A',A\}_-+n\{A',A\}_-\{B',B\}_-=0.
\end{equation}

Remark that the condition $\{A',A\}_-=0$, which may be written as
well as $\frac{A'}A= \frac{\check A'}{\check A}$, is equivalent to
the fact that $\check A= A$. If $a_1$ is $0$, then either
$\{A',A\}_-=0$, which means that $Q=\check P$, or $2\bar b_1
B\check B+n\{B',B\}_-=0$. This last identity is only possible when
$b_1=0$, and thus $\{B',B\}_-=0$. So $Q=P$. In particular, we have
proved the following. At this point, $\pp$ is not necessarily
generic.
\begin{proposition}
Assume that the polynomial $\pp$ is such that $p_1=0$.  Let $\qq$
be a partner of $\pp$. Then $\qq$ is a trivial partner, that is,
there exists a unimodular constant $c$ such that either $\qq=c\pp$
or $\qq=c\check\pp$. \label{p_1=0}\end{proposition}

We will now concentrate on the case when $a_1$ and $b_1$ are
different from zero, and $\pp$ (thus $\qq$) is generic.
 As  $A$ (resp. $B$) has no multiple or symmetric
zeros, then $A \check A$ and $\{A',A\}_-$ (resp. $B \check B$ and
$\{B',B\}_-$) are mutually prime. Moreover, zeros of $A \check A$
and $B \check B$ are different. It follows from (\ref{concl}) that
$2\bar b_1 B\check B+n\{B',B\}_-$ can be divided by $A \check A$,
while $2\bar a_1A \check A+n\{A',A\}_-$ can be divided by $B
\check B$. So $A$ and $B$ have the same degree. We conclude
directly that there is a contradiction when $n$ is odd. From now
on, we  assume that $n=2k$. Then $A$ and $B$ have degree $k$.
Moreover, looking at terms of higher degree, we conclude that
\begin{equation} \label{ab0}
n\{B',B\}_-=2\bar b_1(A\check A -B\check B).\end{equation}
 Differentiating (\ref{ab0}), we obtain
\begin{equation}
\label{ab1} 2\bar b_1(\{A',A\}_+-\{B',B\}_+)=n\{B'',B\}_-.
\end{equation}
We can exchange the roles of $A$ and $B$ in the previous
identities. In particular, we get that
\begin{equation}
\label{ab2} \bar b_1 \{A',A\}_-+\bar a_1\{B',B\}_-=0.
\end{equation}
\end{proof}

\begin{proof}[Second step.]
We will now work with polynomials in two variables. For $\Pi\in
\mathbb{C}[Z,W]$, we define $\check \Pi$ as before, the degree of
a polynomial being taken as the total degree. Using the fact that
$A_\pp \check{A_\pp}=A_\qq \check{A_\qq}$, we know that there
exists a factorization with polynomials $C,D$ in two variables,
such that
\begin{equation}\label{product}
  A_\pp= CD\hspace{2cm} A_\qq= C\check D.
\end{equation}

\medskip

Let us consider $C$ and $D$ as polynomials in the variable $W$
with coefficients that are polynomials in $Z$, and write
\begin{eqnarray*}
C\equiv C_0W^\alpha \hspace{1cm} \mbox{(modulo polynomials in $W$
of
lower degree)}\\
D\equiv D_0W^\beta \hspace{1cm} \mbox{(modulo polynomials in $W$
of lower degree)}.
\end{eqnarray*}
Then $\pp=C_0 D_0 $, while  $\qq= \varepsilon C_0 \check D_0 $,
with $\varepsilon=(-1)^{\deg D + \deg D_0+\beta}$. The assumption
that $\pp$ is generic implies that there is uniqueness in the
factorization (\ref{fact}). So $C_0$ is equal to $A$ (up to a
constant) and $D_0$ is equal to $B$ (up to a constant). Exchanging
the role of the two variables, we see that $\alpha=\beta =k$. So
$\varepsilon=1$, and we can assume that $C_0$ and $D_0$ are monic,
so that $C_0=A$ and $D_0=B$.

 These considerations  allow us to write
\begin{equation}
C(z,w)\equiv A(z)A^*(w)+C_1(z) w^{k-1} \hspace{2cm} D\equiv
B(z)B^*(w)+D_1(z)w^{k-1} \hspace{1cm} \label{first}
\end{equation}
(modulo polynomials in $W$ of lower degree). Moreover, $C_1$ and
$D_1$ have degree at most $k-1$.
\medskip

We shall now identify $A_1$ and $B_1$.

Writing $A_\pp$ as a product, we have that
$$
A_\pp(z,w)\equiv\pp(z)\pp^*(w)+[A(z)D_1(z)+C_1(z)B(z)]w^{n-1}
$$
(modulo polynomials in $W$ of lower degree), whereas a direct
computation, using the fact that $\pp=AB$ shows that
$$
A_\pp(z, w) \equiv \pp(z)\pp^*(w)+n[A(z)B'(z)+A'(z)B(z)]w^{n-1}
$$
(modulo polynomials in $W$ of lower degree). Comparing both
expressions leads to
\begin{equation}\label{C1}
(nA'-C_1)B+(nB'-D_1)A=0.
\end{equation}
Our assumption on the zeros of $\pp$ implies that $A$ and $B$ are
mutually prime so that, using the information on the degrees of
$C_1,D_1$, we get that
$$
C_1=nA'\quad\mbox{and}\quad D_1=nB'.
$$
Symmetry considerations now imply that
\begin{eqnarray*}
C(z,w)\equiv A(z)A^*(w)+2 A^{'}(z)A^{'*}(w)+C_2(z) w^{k-2} \\
D(z,w)\equiv B(z)B^*(w)+2 B^{'}(z)B^{'*}(w)+ D_2(z) w^{k-2}
\end{eqnarray*}
(modulo polynomials in $W$ of lower degree). Moreover, $C_2$
 and $D_2$ have degree at most $k-2$.

\medskip

It then follows that
\begin{eqnarray*}
A_\pp(z,z)&\equiv&\pp(z)\pp^*(w)+\pp'(w){\pp'}^*(w)+\Bigl(A(z)\Bigl[(\bar
a_1-\bar b_1)B'(z)+D_2(z)\Bigl]\\&&
\qquad\qquad\qquad+B(z)\Bigl[(\bar b_1-\bar
a_1)A'(z)+C_2(z)\Bigl]+n^2A'(z)B'(z)\Bigr)w^{n-2}\\
&\equiv&\pp(z)\pp^*(w)+\pp'(z){\pp'}^*(w)
+\frac{n(n-1)}{2}\Bigl(A''(z)B(z)+2A'(z)B'(z)+A(z)B''(z)\Bigr)w^{n-2}
\end{eqnarray*}
(modulo polynomials in $W$ of lower degree). It follows that
\begin{equation}
(\bar a_1-\bar b_1)(AB'-A'B) +AF+BE +nA'B'=0, \label{eq:231}
\end{equation}
where $E:=C_2-\dst\frac{n(n-1)}{2}A''$ and $F:=D_2
-\dst\frac{n(n-1)}{2}B''$. Exploiting the expressions of $A_\qq$,
that is changing $B$ into $\check B$ (thus also $b_1$ into $-b_1$
and $F$ into $\check F$), we get
\begin{equation}
(\bar a_1+\bar b_1)(A\check B'-A'\check B) +A\check F+\check B E
+nA'\check B'=0. \label{eq:232}
\end{equation}
Let us multiply the left hand side of (\ref{eq:231}) by $\check B$
and the left hand side of (\ref{eq:232}) by $B$, and take the
difference. We obtain that $$A\left(\bar a_1 \{B',B\}_--\bar b_1
\{B',B\}_++\{F,B\}_- \right)+A'\left (2\bar b_1 B\check B+n
\{B',B\}_+\right)=0.$$ Using (\ref{ab0}) and (\ref{ab2}), we can
write that
$$A'\left (2\bar b_1 B\check B+n
\{B',B\}_+\right)=2\bar b_1 AA'\check A= A\left(\bar b_1
\{A',A\}_+-\bar a_1 \{B',B\}_-\right ).$$ Finally, using
(\ref{ab1}), we obtain the identity
$$\{F+ \frac n2 B'',B\}_- =0.$$
Since $B$ and $\check B$ are mutually prime by assumption, this
means that $F=-\frac n2 B''$.

We could as well prove that $E=-\frac n2 A''$. If we compute the
coefficient of the term of higher degree in the left hand side of
(\ref{eq:231}), we obtain $ |a_1|^2+ |b_1|^2+n$, which cannot
vanish.
 This concludes for the proof.
\end{proof}
\begin{remarknum}
\label{rem:generic} We will need the following: for all integers
$n\neq m$ and $a\in\mathbb{C}$, then $\psi_n+a\psi_m$ has only
trivial partners for the ambiguity problem (equivalently,
$Z^n+aZ^m$ has only trivial partners for the algebraic ambiguity
problem). Indeed, for $|n-m|\geq 2$ this is a consequence of
Proposition \ref{p_1=0}. For $|n-m|=1$, this follows directly from
(\ref{P=Q}).
\end{remarknum}

\section{Trivial solutions and constructions of special strange partners}
\label{sec:3}

\subsection{The discrete case}
We refer to Problem 3 as Problem (P). In this setting, two
sequences $a$ and $b$ are said to be {\it discrete ambiguity
partners (or}  (P)-{\it partners)} whenever (\ref{pbdisc}) holds.

We start by defining the dual problem of (P), when $2\pi$-periodic
functions, rather than sequences in $\Z$, are considered. Here
$\T=\R/2\pi\Z\equiv[0,2\pi)$, and for $f\in L^2(\T)$ we let
$$
\hat f(n)=\frac1{2\pi}\int_0^{2\pi}f(t)\,e^{-int}\,dt,\quad
n\in\Z.
$$
In this way, one can write $f(t)=\sum_{n\in\Z}\hat f(n)\,e^{int}$
in the usual $L^2(\T)$ sense (and $a.e.$). We shall also identify
$L^2(\T)$ with $\ell^2(\Z)$ via the correspondence: $f\mapsto
\{\hat f(n)\}_{n\in\Z}.$ This gives the following equivalent
formulation of (P).

\medskip
\noindent{\bf ($\mathbf{\widehat{\bf P}}$)\ The Periodic Ambiguity
Problem.} {\sl For $f\in L^2(\T)$ define the periodic  ambiguity
function by
$$
\widehat{\aa}(f)(k,t)=\frac1{2\pi}\int_0^{2\pi}f(s){\Ol{f(s-t)}}\,e^{-iks}\,ds
,\quad (k,t)\in\Z\times\T.
$$
We want to find all $g\in L^2(\T)$ such that
$$
\abs{\widehat{\aa}(f)(k,t)}=\abs{\widehat{\aa}(g)(k,t)},\qquad\mathrm{for\
all\ }(k,t)\in\Z\times\T.
$$
Two functions $f$ and $g$ as above are called {\rm
$\hP$-partners}.}

Note that $f$ and $g$ are  $\hat P$-partners if and only if the
sequences of their Fourier coefficients $\{\hat f(n)\}$ and
$\{\hat g(n)\}$ are (P)-partners in the sense of (\ref{pbdisc})
since Parseval's formula gives
\begin{equation}
\label{parseval}
  \widehat{\aa}(f)(k,t)=\sum_{n\in\Z}\hat f(n){\overline{\hat f(n-k)}}e^{int}
=\aa\bigl(\{\hf(n)\}\bigr)(k,t).
\end{equation}
In the sequel, we will therefore write $\aa(f)$ instead of
$\widehat{\aa}(f)$ to simplify the notation.

\medskip
Let us give a precise definition of  trivial solutions, as
announced in the introduction. Intuitively these should be  simple
transformations of the data function that always give solutions to
the functional equation proposed. The definition below, given for
$\hP$ easily adapts to other problems.

\begin{defi} A \emph{trivial solution for $\hP$} is a
 bounded linear operator $R\,:\ L^2(\T)\to L^2(\T)$
preserving $\hP$-partners, {\it i.e.}, such that for every $f\in
L^2(\T)$, $f$ and $Rf$ are $\hP$-partners. We denote by
$\mathcal{T}$ the semi-group of all such operators.
\end{defi}

\begin{example} Let $\hil\equiv\T\times\Z\times\T$, and define for
$h=(\al,k,\beta)\in\hil$
$$
R_hf(t)=e^{i\beta}e^{ikt}f(t+\al)\quad\mbox{and}\quad
\tilde R_hf(t)=e^{i\beta}e^{-ikt}f(-t+\al).
$$

Then $R_h$ and $\tilde R_h$ are trivial solutions for $\hP$. Note
that $R_h$ is a unitary representation of the periodized
Heisenberg group $\hil$, with the product defined by
$$
h\cdot
h'=(\al,k,\beta)\cdot(\al',k',\beta')=(\al+\al',k+k',\beta+\beta'+k'\al)
$$
while as before $\tilde R_h=ZR_h$.
\end{example}

Let us prove that there are no other trivial solutions.

\begin{prop}
\label{tl2} Let $R$ be a trivial solution for $\hP$, then there
exists $h\in\hil$ such that, either $R=R_h$ or $R=\tilde R_h$. In
particular, $\mathcal{T}$ can be identified with the group
$\{-1,1\}\times\hil$.
\end{prop}

\begin{proof}[Proof]
For $n\in\Z$, let $f_n(t)=e^{int}$. Then,
$\abs{\aa(f_n)(k,t)}=\dt_{0,k}$, where $\dt_{0,k}$ is the usual
Kr\"onecker symbol. Moreover,
$$
\abs{\aa(Rf_n)(k,t)}
=\abs{\sum_{\ell\in\Z}\widehat{Rf_n}(\ell)\overline{\widehat{Rf_n}(\ell-k)}e^{i\ell
t}}=\dt_{0,k}
$$
implies that there exists a unique $m(n)$ such that
$\widehat{Rf_n}\bigl(m(n)\bigr)\not=0$, that is
$Rf_n(t)=c_ne^{im(n)t}$ with $|c_n|=1$. Note that, if
$n_1\not=n_2$, then  $m(n_1)$ and $m(n_2)$ are different. Indeed,
if they were equal, the non-zero function
$g:=c_{n_2}f_{n_1}-c_{n_1}f_{n_2}$ would have a zero radar
ambiguity function, a clear contradiction.

We wish to show that either $m(n)-n$ or $m(n)+n$ is a constant.
Let us consider the test functions $g(t)=e^{in_1t}+e^{in_2t}$, for
distinct $n_1,n_2\in\Z$. Then
$Rg(t)=c_{n_1}e^{im(n_1)t}+c_{n_2}e^{im(n_2)t}$, and therefore,
$$
\abs{\aa(g)(0,t)}=
\abs{e^{in_1t}+e^{in_2t}}=\abs{\abs{c_{n_1}}^2e^{im(n_1)t}
+\abs{c_{n_2}}^2e^{im(n_2)t}}= \abs{\aa(Rg)(0,t)}.
$$

This implies, $|m(n_1)-m(n_2)|=|n_1-n_2|$, which is an isometry of
the integers, and therefore of the form $m(n)=m(0)+\eps n$, with a
constant $\eps=\pm1$. In particular, when $\eps=1$ we have
\begin{equation}
(Rf_n)(t)=c_n e^{im(0)t} f_n(t). \label{mul}
\end{equation}
We shall show that actually $R=R_h$ for some $h\in\hil$. The case
$\eps=-1$, then follows by replacing $R$ by $R Z$.

So, assuming (\ref{mul}), let us establish the dependence of $c_n$
on $n$. Testing with $h_n(t)=e^{int}+e^{i(n+1)t}+e^{i(n+2)t}$, we
obtain
\begin{align}
|\aa(h_n)(1,t)|=  |1+e^{it}|= |1+c_n\overline{c_{n+1}}^2
c_{n+2}e^{it}|=|\aa(Rh_n)(1,t)| .\notag
\end{align}
 Therefore $\overline{c_{n+1}}c_{n+2}=\overline{c_n}c_{n+1}$. Writing
$c_n=e^{i\gamma(n)}$, this relation can be expressed as
$$
\gamma(n+2)-\gamma(n+1)=\gamma(n+1)-\gamma(n)=\ldots=
\gamma(1)-\gamma(0)\qquad (mod\ 2\pi).
$$
Hence, for some $\al\in\T$, we must have
$$
\gamma(n)=\al+\gamma(n-1)=\ldots=n\al+\gamma(0)\qquad (mod\ 2\pi),
$$
concluding that
$$
(Rf_n)(t)=c_0e^{in\al}e^{im(0)t}f_n(t)=c_0e^{im(0)t}f_n(t+\al).
$$
Then, the linearity and boundedness of $R$  give $R=R_h$, where
$h=(\al,m(0),\gamma(0))$.
\end{proof}

\begin{remarknum}
It is worthwhile to notice that, from the above proof, an
anti-linear bounded operator $R$ cannot preserve $\hP$-partners.
Indeed, in the last step of the proof one may test with a function
$f(t)=1+e^{it}+ce^{2it}$, for $|c|=1$. Then
$\abs{\aa(f)(1,t)}=\abs{1+c e^{it}}$, whereas if $R$ was
antilinear, $\abs{\widehat{\aa}(f)(1,t)}
=\abs{\aa(Rf)(1,t)}=\abs{1+\bar ce^{it}}$. This excludes
anti-linear operators to give trivial solutions for $\hP$.
\end{remarknum}

\noindent{\bf A normalization remark.} \,
Let $f\in L^2(\T)$ be a trigonometric polynomial. Then, up to a change
$f\mapsto e^{ikt}f$, we may assume that $\supp\widehat{f}\subset\{0,\ldots,N\}$ for some integer $N$
and that $\widehat{f}(0)\not=0$, $\widehat{f}(N)\not=0$. We then
say that $f\in\pp_N$.

The next lemma shows in particular that there is no loss of
generality if we restrict the study of the discrete radar
ambiguity problem to functions in $\pp_N$ when dealing with
trigonometric polynomials.

\begin{lemma}
\label{tl3} Let $f\in L^2(\T)$ and let $\Lambda=\supp f$. Then
$$
\supp\aa(f):=\{k\tq\aa(f)(k,t)\ \mathrm{is\ not\ identically\
}0\}=\Lambda-\Lambda.
$$
In particular, if $f\in\pp_N$ for some $N\in\N$, and if $g$ is a $(\widehat{P})$-partner of $f$
then, up to replacing $g$ by a trivial partner, we may also assume that $g\in\pp_N$.
\end{lemma}

\begin{proof} The $n$-th Fourier coefficients of $t\mapsto\aa(f)(k,t)$, namely
$\hf(n){\Ol{\hf(n-k)}}$, will vanish unless $n,n-k\in\La$, so that
$\supp\aa(f)=\Lambda-\Lambda$.

If $f\in\pp_N$ then $\Lambda\subset\{0,\ldots,N\}$, thus
$\supp\aa(f)\subset\{0,\ldots,N\}-\{0,\ldots,N\}=\{-N,\ldots,N\}$.
Obviously $\aa(f)(-N,t)=\widehat{f}(0)\overline{\widehat{f}(N)}\not=0$,
$\aa(f)(N,t)=\widehat{f}(N)\overline{\widehat{f}(0)}e^{iNt}\not=0$,
thus $\supp\aa(f)$ cannot be included in a smaller interval.

Now, if $g\in L^2(\T)$ is a $(\widehat{P})$-partner of $g$, then
$\Lambda':=\supp\widehat{g}$ is such that $\Lambda'-\Lambda'$ is finite, thus $\Lambda'$ itself
is finite. Thus $g$ is a trigonometric polynomial, thus we may assume that $g\in\pp_M$ for some $M$.
The first part of the proof then shows that $M=N$.
\end{proof}

Finally, it is obvious from the definition that if $f\in\pp_N$ with
$N=0$ or $N=1$, then $f$ has only trivial partners.

\subsection{Restricted discrete problems}

In this section we consider the discrete radar ambiguity problem
($\hat{P}$) restricted to the subspaces $\Ll$. Recall that, for
$\La$ a subset of $\mathbb{Z}$, this space consists of all
functions $f\in L^2(\T)$ with $\supp\hat f\subset\La$. The
discrete radar ambiguity problem may then be restricted in two
ways:

\medskip
\noindent{\bf The Ambiguity Problems in $\Ll$.} {\sl Given
$f\in\Ll$, }
\begin{itemize}
\item[$\mathbf{\widehat{P}_\Lambda}$.] {\sl find all $g\in\LT$ such that
for all $(k,t)\in\Z\times\T$}
\begin{equation}
\label{ramba} \tag{$\widehat{P}_\Lambda$} \abs{\aa
(f)(k,t)}=\abs{\aa (g)(k,t)}
\end{equation}
and such a $g$ will be called a $\widehat{P}_\Lambda$-partner of
$f$;

\item[$\mathbf{\widehat{P}_{\Lambda,\Lambda}}$.] {\sl  find all $g\in\Ll$ such that
for all $(k,t)\in\Z\times\T$}
\begin{equation}
\label{rambb} \tag{$\widehat{P}_{\Lambda,\Lambda}$} \abs{\aa
(f)(k,t)}=\abs{\aa (g)(k,t)}.
\end{equation}
Such a $g$ will be called a
$\widehat{P}_{\Lambda,\Lambda}$-partner of $f$.
\end{itemize}
In other words, the $\widehat{P}_\Lambda$-ambiguity problem is
just the $\hat{P}$-ambiguity problem for functions in $\Ll$
whereas in the $\widehat{P}_{\Lambda,\Lambda}$-ambiguity problem
one further seeks for the solutions of the $\hat{P}$-ambiguity
partners to be in $\Ll$.
\medskip
Restricted trivial solutions may now be defined in two natural
ways:
\begin{itemize}
\item an operator $R\,:\ \Ll\to\LT$ such that, for every $f\in\Ll$,
$f$ and $g=Rf$ satisfy (\ref{ramba}) will be called a {\it
$\widehat{P}_\Lambda$-trivial solution};

\item an operator $R\,:\ \Ll\to\Ll$ such that, for every $f\in\Ll$,
$f$ and $g=Rf$ satisfy (\ref{rambb}) will be called a {\it
$\widehat{P}_{\Lambda,\Lambda}$-trivial solution}.
\end{itemize}

Of course, every $\widehat{P}_{\Lambda,\Lambda}$-trivial solution
is also a $\widehat{P}_\Lambda$-trivial solution. The converse may
not be true as the trivial solutions $R_{0,k,0}$ and
$\widetilde{R}_{0,k,0}$ do not preserve $\Ll$ in general. Note
also that every trivial solution is a
$\widehat{P}_\Lambda$-trivial solution. Again the converse may be
false as the example bellow will show.

It is a  remarkable fact that the more {\it lacunary} a sequence
$\La$ is, the more trivial solutions the problem admits.

\begin{notation} For $\Lambda\subset\Z$ and $\bc=\{c(n)\}_{n\in\La}$
a sequence of unimodular numbers, we define the (multiplier)
operator $R_\bc:\Ll\to\Ll$ by
$$
R_\bc f(t)=\sum_{n\in\La}c(n)\hat{f}(n)e^{int},\quad t\in\T.
$$
This operator is extended to $L^2(\T)$ in the obvious way: $R_\bc
e^{int}=0$ if $n\notin\Lambda$.
\end{notation}

\begin{example}
Let $\La=\{2^j\}_{j=0}^\infty$. Then, any multiplier $R_\bc$ is a
trivial solution for $\Pl$, but in general not for $(\hat{P})$.
This is due to the fact that $\aa f(0,t)=\aa R_\bc f(0,t)$, while
$$
|\aa f(2^{k_1}-2^{k_2},t)|
=\abs{\hf(2^{k_1}){\Ol{\hf(2^{k_2})}}}=|\aa R_\bc
f(2^{k_1}-2^{k_2},t)|,
$$
for non-negative integers $k_1\not= k_2$.
\end{example}

In general, we have the following result:

\begin{prop}
\label{prop1} Let $\La\subset\Z$. An operator $R\,:\ \Ll\to\LT$ is
a $\Pl$-trivial solution if and only if it is of the form
$R=SR_\bc$, where $S\in\mathcal{T}$ and $\bc=\{c(n)\}_{n\in\La}$
is a sequence of unimodular constants satisfying
\begin{equation}
\label{con} c(n_1){\Ol{c(n_2)}}=c(n_3){\Ol{c(n_4)}},\quad
\mbox{whenever}\quad n_1-n_2=n_3-n_4,\,\, n_i\in\La,\,\,\,
i=1,2,3,4.
\end{equation}
\end{prop}

\begin{proof}
The sufficiency is easy to check. Indeed, just notice that for
$R=R_{\bc}$, using condition (\ref{con}) one obtains
$$
\abs{\aa(Rf)(k,t)}  = \abs{\sum_{n,n-k\in\La}
c(n)\hf(n){\Ol{c(n-k)\hf(n-k)}}e^{int}}
 =\abs{\aa(f)(k,t)}
$$
since $c(n)\Ol{c(n-k)}$ depends only $k$ and is of modulus $1$.

For the necessity, it is easy to see that the operator $R$ will
act on the exponentials $f_n(t)=e^{int}$ by
$$
\mbox{either}\quad Rf_n(t)=c(n)e^{imt}f_n(t),\quad \mbox{or}\quad
Rf_n(t)=c(n)e^{-im t}f_{-n}(t),
$$
for some $m\in\Z$ and $|c(n)|=1$, $n\in\La$. Indeed, the part of
the proof of Proposition \ref{tl2} to give (\ref{mul}) can be used
here. Factoring out the corresponding $\hP$-trivial operator
$S=R_{h}$ or $S=\tilde{R}_h$ with $h=(0,m,0)\in\hil$, we may
assume that $R=R_{\bc}$. It remains to determine the relations in
(\ref{con}) among the $c(n)$'s.

Excluding the trivial cases, we have only to check (\ref{con})
when $n_1>n_2$ and $n_1>n_3>n_4$. This leaves only two
possibilities:

\medskip
\noindent{\bf Case 1:} $n_2=n_3$. Then testing with
$g(t)=e^{in_1t}+e^{in_2t}+e^{in_4t}$, we obtain
$$
\abs{\aa(Rg)(n_1-n_2,t)} = \abs{c(n_1){\Ol{c(n_2)}}e^{in_1t}
+c(n_2){\Ol{c(n_4)}}e^{in_2t} }
 =\abs{e^{in_1t}+e^{in_2t}},
$$
and consequently, $c(n_1){\Ol{c(n_2)}}=c(n_2){\Ol{c(n_4)}}$.

\medskip
\noindent{\bf Case 2:} $n_2\neq n_3$. Then, the $n_i$'s are all
different and we may test
 with $h(t)=e^{in_1t}+e^{in_2t}+e^{in_3t}+e^{in_4t}$,
obtaining:
$$
\abs{\aa(Rh)(n_1-n_3,t)}
 =\abs{e^{in_1t}+e^{in_2t}+\th(n_3){\Ol{\th(2n_3-n_1)}}e^{in_3t}}.
$$
Note that $\th(2n_3-n_1)\not=0$ only if $2n_3-n_1=n_2$ or  $n_4$.
But the last choice implies $n_2=n_3$, which is not possible. If
instead $n_3-n_1=n_2-n_3$, then the previous case gives us the
equality
$$
c(n_3){\Ol{c(n_1)}}=c(n_2){\Ol{c(n_3)}}.
$$
Therefore
$$
\abs{\aa(Rh)(n_1-n_3,t)} =
\abs{c(n_1){\Ol{c(n_3)}}(e^{in_1t}+e^{in_3t})
+c(n_2){\Ol{c(n_4)}}e^{in_2t} }
 =\abs{e^{in_1t}+e^{in_2t}+e^{in_3t}},
$$
from which we obtain $c(n_1){\Ol{c(n_3)}}=c(n_2){\Ol{c(n_4)}}$.
When, on the contrary, $\th(2n_3-n_1)=0$, then the situation is
simpler since
$$
\abs{\aa(Rh)(n_1-n_3,t)} = \abs{c(n_1){\Ol{c(n_3)}}e^{in_1t}
+c(n_2){\Ol{c(n_4)}}e^{in_2t} }
 =\abs{e^{in_1t}+e^{in_2t}},
$$
leading to the same result.
\end{proof}

\begin{remark}
We shall denote by $\mathcal{T}_\La$ the set of all operators
which are trivial solutions for $\Pl$. Note that now
$\mathcal{T}_\La$ is not a semi-group with the usual composition
law, unless $\La=\Z$.
\end{remark}

The previous proposition, translated into the language of the
periodic  radar ambiguity problem $\hP$, guarantees the existence
of {\it many strange solutions} for every function in $\Ll$,
 provided $\La$ has
enough gaps. Further, we obtain the following:

\begin{cor}
The set of functions $f\in L^2(\T)$ admitting strange solutions to
$\hP$ is dense in $L^2(\T)$.
\end{cor}

\begin{proof}
Consider, for every $N\geq1$, functions with Fourier transform
supported in $\La_N=\{-N,\ldots,N\}\cup\{3N+1\}$. It is clear that
$\cup_{N=1}^\infty L^2_{\La_N}(\T)$ is dense in $L^2(\T)$.
Further, any function $f\in L^2_{\La_N}(\T)$ will have infinitely
many $\hP$-strange partners. Indeed, these are given by the
$\Pl$-trivial solutions:
$$
R_{\bc_N}f(t)=\sum_{n=-N}^N\eps \hf(n)e^{int}+
\eps'\hf(3N+1)e^{i(3N+1)t},
$$
for $|\eps|=|\eps'|=1$. Since the multiplier
$\bc_N=\{c(-N)=\ldots=c(N)=\eps,c(3N+1)=\eps'\}$ satisfies
condition (iii) of Proposition \ref{prop1}, we must have
$R_{\bc_N}\in\mathcal{T}_{\La_N}$, establishing our claim.
\end{proof}

Here is one more consequence of our proposition, generalizing the
example given above. We exclude the case Card$(\La)=2$ for which
one easily knows all solutions to $\hP$ or $\Pl$.

\begin{cor}
\label{corlac} Let $\La\subset\Z$ be such that Card$(\La)\geq3$.
Suppose that every $n\in\La+\La$ can be written uniquely (up to
permutation) as $n=n_1+n_2$, with $n_1,n_2\in\La$. Then
$R:\Ll\to\Ll$ is a $\widehat{P}_{\Lambda,\Lambda}$-trivial
solution if and only if it is of the form $R=R_{\bc}$ with
$\bc\equiv\{c(n)\}_{n\in\La}\in\T^\La$.

Further, if $f\in\Ll$ and Card$(\supp\hf)\geq3$, then every
solution to $\widehat{P}_{\Lambda,\Lambda}$ is given by $R_\bc f$,
for some $\bc\in\T^\La$.
\end{cor}

In other words, this corollary states that the trivial solutions
may be identified with $\T^\Lambda$ and that, if $f\in\Ll$ and
Card$(\supp\hf)\geq3$, every solution to
$\widehat{P}_{\Lambda,\Lambda}$ is a
$\widehat{P}_{\Lambda,\Lambda}$-trivial solution.

\begin{proof}
Under the assumption on $\La$, condition (iii) of Proposition
\ref{prop1} always holds, since $n_1-n_2=n_3-n_4$, for $n_i\in\La$
implies $n_1=n_3$ or $n_1=n_2$. It follows that the $\Pl$-trivial
solutions are all given by $R_\bc R_{\alpha,k,\beta}$ or by
$R_\bc\widetilde{R}_{\alpha,k,\beta}$ for some
$(\alpha,k,\beta)\in\mathbb{H}$ and some
$\bc=\{c(n)\}_{n\in\Lambda}\in(S^1)^\La\equiv\T^\La$. Among these
operators, the only ones that preserve $\Ll$ are
$R_\bc\widetilde{R}_{\alpha,0,\beta}=R_{\tilde\bc}$ with
$\tilde\bc_n=e^{i\beta+in\alpha}\bc_n$.

We shall show that if $g\in\Ll$ is a
$\widehat{P}_{\Lambda,\Lambda}$-partner of $f$, and
Card$(\supp\hf)\geq3$, then $|\hf(n)|=|\hg(n)|$, for all
$n\in\La$. This will imply that $g=R_\bc f$ for some multiplier
$\bc\in\T^\La$ and establish the corollary.

From the assumptions on $\La$,  it follows that
$|\aa(f)(k,t)|=|\aa(g)(k,t)|$ is a constant for each $k\in\Z$. For
instance, if we fix $n_0\in\supp\hf$, then for every
$n\in\La\setminus\{n_0\}$, we get, for $k=n_0-n\not=0$
\begin{equation}
\label{equal}
\bigl|\hf(n_0){\Ol{\hf(n)}}\bigr|=\bigl|\hg(n_0){\Ol{\hg(n)}}\bigr|.
\end{equation}
Since Card$(\supp\hf)\geq2$, we must have $\hg(n_0)\not=0$.
Denoting $z=\frac{\hf(n_0)}{\hg(n_0)}$, and using
$|\aa(f)(0,t)|^2=|\aa(g)(0,t)|^2$ we obtain:
$$
\abs{\abs{\hf(n_0)}^2e^{in_0t}
+\sum_{n\not=n_0}\abs{\hf(n)}^2e^{int}}^2=
\abs{\frac{1}{\abs{z}^2}\abs{\hf(n_0)}^2e^{in_0t}
+\abs{z}^2\,\sum_{n\not=n_0}\abs{\hf(n)}^2e^{int}}^2.
$$
Thus both trigonometric polynomials have same coefficients so that
either $|z|=1$ (which is what we wish), or
$$
\abs{\hf(\la_0)}^4=
|z|^4\abs{\sum_{\la\not=\la_0}|\hf(\la)|^2e^{i\la t}}.
$$
Since Card$(\supp\hf)\geq3$, we see that the latter cannot happen,
so that $|z|=1$. The corollary then follows from (\ref{equal}).
\end{proof}
\begin{remarknum}\label{two} Sets $\Lambda$ satisfying the condition of the corollary
are usually called $B_2$-sets (or $B_2[1]$-sets) and have been
extensively studied by Erd\"os and various collaborators, as well
as their generalization, the $B_k$-sets, where sums of two
integers are replaced by sums of $k$ integers. One may show that a
subset $\Lambda$ of $\{1,\ldots,N\}$ that is a $B_k$ set has size
at most $\mbox{Card}\,\Lambda\leq CN^{1/k}$ and this bound is
sharp. A survey on the subject may be found on M. Koluntzakis' web
page (see also \cite{Ko1,Ko2}). $B_k$-sets are particular examples
of {\it $\La(2k)$-sets} for trigonometric series \cite{rudin}. The
two dimensional version of these sets, contained in the lattice
$\Z^2$, also appears in the study of certain phase retrieval
problems arising from crystallography \cite{tes}.
\end{remarknum}

When the gaps of $\La$ are even larger, we will now prove that the
problem $\widehat{P}_\Lambda$ has {\it only} trivial solutions. To
do so, we will need the following lemma which may be well known.

\begin{lemma}
\label{lem:b3} Let $\Lambda,\Lambda'\subset\Z$ and assume that
every $n\in\Lambda+\Lambda+\Lambda$ can be written uniquely up to
permutation as $n=n_1+n_2+n_3$ with $n_1$, $n_2$, $n_3\in\Lambda$.
Assume further that $\Lambda'-\Lambda'=\Lambda-\Lambda$, then
$\Lambda'=\Lambda-m$ or $\Lambda'=m-\Lambda$ for some $m\in\Z$.
\end{lemma}

\begin{proof} Without loss of generality, we may assume that $0\in\Lambda,\Lambda'$. Now, if
$m\in\Lambda'\setminus\{0\}$, we may write $m=m-0=n_1-n_2$ for
some $n_1,n_2\in\Lambda$. Assume that we may write
$m=n_1^{'}-n_2^{'}$ with $n_1^{'},n_2^{'}\in\Lambda$, then
$n_1+n_2^{'}+0=n_1^{'}+n_2+0$. The property of $\Lambda$ together
with $m\not=0$ then implies that $n_1^{'}=n_1$ and $n_2^{'}=n_2$.
It follows that every $m\in\Lambda'\setminus\{0\}$ may be written
in a unique way as $m=n_m-\widetilde{n}_m$ with
$n_m\not=\widetilde{n}_m\in\Lambda$.

Further, fix $m_0\in\Lambda'\setminus\{0\}$ and write
$m_0=n_0-\widetilde{n}_0$ with
$n_0\not=\widetilde{n}_0\in\Lambda$. Then, for
$m\in\Lambda'\setminus\{0,m_0\}$, as
$m-m_0\in\Lambda'-\Lambda'=\Lambda-\Lambda$, there exist
$n\not=\widetilde{n}\in\Lambda$ such that $m-m_0=n-\widetilde{n}$.
It follows that
$n_m+\widetilde{n}_0+\widetilde{n}=\widetilde{n}_m+n_0+n$. As
$m\not=0$, we get $\widetilde{n}_m\not=n_m$ and as $m\not=m_0$, we
get $\widetilde{n}\not=n$. The condition on $\Lambda$ then implies
that either
$(n_m,\widetilde{n}_0,\widetilde{n})=(n_0,n,\widetilde{n}_m)$ or
$(n_m,\widetilde{n}_0,\widetilde{n})=(n,\widetilde{n}_m,n)$. In
the first case, $m=n_m-\widetilde{n}_m=n_0-\widetilde{n}_m\in
n_0-\Lambda$ while in the second case
$m=n_m-\widetilde{n}_m=n_m-\widetilde{n}_0\in\Lambda-\widetilde{n}_0$.

It is now enough to prove that, for a given $\Lambda'$, only one
of these cases may occur.

If $\mbox{Card}\,\Lambda'\leq2$ this is trivial. If
$\mbox{Card}\,\Lambda'=3$, the uniqueness of the decomposition
$0\not=m=n_m-\widetilde{n}_m$ implies that, if
$m\in(\Lambda-\widetilde{n}_0)\cap(n_0-\Lambda)$ then $m=m_0$. We
may thus assume that $\mbox{Card}\,\Lambda'\geq4$.

Let $m\not=\widetilde{m}\in\Lambda'\setminus\{0,m_0\}$ and assume
that we may write $m=n_0-n$ and
$\widetilde{m}=\widetilde{n}-\widetilde{n}_0$ with
$n,\widetilde{n}\in\Lambda$. Again, as
$\Lambda'-\Lambda'=\Lambda-\Lambda$, there exists $n_1\not=n_2$
such that $m-\widetilde{m}=n_1-n_2$. It follows that
$n_0+\widetilde{n}_0+n_2=n+\widetilde{n}+n_1$. The property of
$\Lambda$ with $n_1\not=n_2$ then implies that only four cases may
occur:
$$
(n_0,\widetilde{n}_0,n_2)=(n,n_1,\widetilde{n}),\
(n_0,\widetilde{n}_0,n_2)=(n_1,\widetilde{n},n),\
(n_0,\widetilde{n}_0,n_2)=(\widetilde{n},n_1,n)\mbox{ or }
(n_0,\widetilde{n}_0,n_2)=(n_1,\widetilde{n},n).
$$
The two first cases are respectively excluded with $m\not=0$ {\it
i.e.} $n_0\not=n$ and $\widetilde{m}\not=0$ {\it i.e}
$\widetilde{n}_0\not=\widetilde{n}$. The two last cases are
respectively excluded with $\widetilde{m}\not=m_0$ {\it i.e.}
$n_0\not=\widetilde{n}$ and $m\not=m_0$ {\it i.e.}
$\widetilde{n}_0\not=n$. This concludes the proof of the lemma.
\end{proof}

 Sets $\Lambda$ satisfying the condition of the lemma
are usually called $B_3$-sets. See Remark \ref{two} above.

\begin{cor}
\label{corlac2} Let $\La\subset\Z$ be a $B_3$-set. Then every
solution to $\widehat{P}_\Lambda$ is a trivial solution, that is
if $f\in\Ll$, then the solutions to $\hP$ are all given by $SR_\bc
f$, for $\bc\in\T^\La$, $S\in\mathcal{T}$.
\end{cor}

\begin{proof}
Without loss of generality we assume $0\in\supp\hf$ and
Card$(\supp\hf)\geq3$. Note that, since $\La$ satisfies the
assumptions in Corollary \ref{corlac},  all the solutions to
$(\widehat{P}_{\Lambda,\Lambda})$ are given by $R_\bc f$.

We shall show that if $g$ is a $\hP$-partner of $f$, then
$\supp\widehat{Sg}\subset\La$, for some $S\in\mathcal{T}$. This
will imply that $f$ and $Sg$ are
$(\widehat{P}_{\Lambda,\Lambda})$-partners, and hence
$g=S^{-1}R_\bc f$.

We denote $\Lf=\supp\hf$ and $\Lg=\supp\hg$. As $f$ and $g$ are
ambiguity partners, $\aa(f)$ and $\aa(g)$ have same support and,
with Lemma \ref{tl3} this implies that $\Lf-\Lf=\Lg-\Lg$. From
Lemma \ref{lem:b3}, we get that either $\Lg=\Lf-m$ or $\Lg=m-\Lf$
for some $m\in\Z$.

In the first case, it suffices to define $S\in\mathcal{T}$ by
$Sg(t)=e^{-imt}g(t)$ while in the second case we consider
$Sg(t)=e^{imt}g(-t)$. We then have
$\supp\widehat{Sg}\subset\supp\hf$ and, hence, $f$ and $Sg$ are
$(\widehat{P}_{\Lambda,\Lambda})$-partners. The proof of the
corollary is then complete.
\end{proof}

To conclude this section, let us point out the existing relation
between $\Pl$-trivial solutions and ``restricted'' solutions to
the ambiguity problem, as they were defined for
the continuous case (\ref{rap}) in \cite{Ja}.
In the periodic situation, the question can be asked as follows:

\medskip
\noindent{\bf ($\mathbf{\widehat{\bf P_r}}$\bf) The Restricted Ambiguity Problem.}
{\sl For $f\in L^2(\T)$,
find all $g\in L^2(\T)$ for which there is some family of unimodular constants $\eta_k$
such that, for all $(k,t)\in\Z\times\T$
\begin{equation}
\widehat{\aa}(f)(k,t)=\eta_k\widehat{\aa}(g)(k,t).
\label{res}
\end{equation}
Two functions $f$ and $g$ as above are  called
{\it restricted partners}.}
\medskip

We have the following result:

\begin{cor}
\label{restr}
Let $f\in L^2(\T)$ and $\La=\supp\hf$. Then, all the restricted partners of $f$
are of the form
$R_{\bc}f$, with $R_\bc\in\mathcal{T}_\La$, that is $\bc$ is a sequence of unimodular
constants supported in $\Lambda$ that satisfies (\ref{con}).
\end{cor}

\begin{proof}
It is clear that for each $R_{\bc}\in\mathcal{T}_\La$, with
$\La=\supp\hf$, then $R_\bc f$ is a restricted partner of $f$.
Indeed, (\ref{res}) holds with $\eta_k=c(n){\Ol{c(n-k)}}$,
which by (\ref{con}) does not depend on $n\in\La$.
Conversely, Equality (\ref{res}) for $k=0$ implies
$|\hf(n)|=|\hat{g}(n)|$ for all $n\in\Z$.
Thus, $g\in\Ll$ and $g=R_\bc f$
for a sequence of unimodular constant $\bc=\{c(n)\}_{n\in\La}$.
It remains to show that condition (iii) in Proposition
\ref{prop1} holds. But this  once more follows from (\ref{res}),
since for general values of $k\in\La-\La$, have
$\eta_k={\Ol{c(n)}}c(n-k)$, for all
$n,n-k\in\La$.
\end{proof}

\subsection{The continuous case}

The definition of trivial solutions immediately adapts to the
continuous radar ambiguity problem: a \emph{trivial solution} to
the continuous radar ambiguity problem is a linear or anti-linear
continuous operator $T$ on $\Lt$ such that for every $u\in\Lt$,
$u$ and $Tu$ are ambiguity partners. We have the following
description of these operators:

\begin{proposition} The trivial solutions of the continuous radar ambiguity
are the operators of the form $Tu(t)=ce^{i\omega t}u\bigl(\eps(t-a)\bigr)$
with $c\in\T$, $\eps=\pm1$, $\omega,a\in\R$.
\end{proposition}

\begin{proof} Let $T$ be a trivial solution and let $\psi_n$ be the Hermite basis.
According to Remark \ref{rem:generic}, $\psi_n$, $\psi_n+\psi_k$ have only trivial partners.
Thus, for every $n$,
there exists  $c_n\in\T$, $\eps_n=\pm1$, $\omega_n,a_n\in\R$ such that
$$
T\psi_n(t)=c_ne^{i\omega_nt}\psi_n\bigl(\eps_n(t-a_n)\bigr)=
c_n\eps_n^ne^{i\omega_nt}\psi_n(t-a_n).
$$
We want to prove that these constants do not depend on $n$:
$a_n=a_0,\omega_n=\omega_0$ and either $c_n\eps_n^n=c_0$ or
$c_n\eps_n^n=(-1)^nc_0$. If this is the case, then respectively
$T\psi_n(t)= c_0e^{i\omega_0t}\psi_n(t-a_0)$ or
$T\psi_n(t)= c_0e^{i\omega_0t}\psi_n(-t+a_0)$. By density of the span of the
$\psi_n$'s, linearity and continuity of $T$,
it follows that $Tu(t)=c_0e^{i\omega_0t}u\bigl(\eps_1(t-a_0)\bigr)$ for all
$u\in L^2$, as desired.

To do so, take $n\not=k$ and note that by additivity of $T$,
$$
T(\psi_n+\psi_k)=T\psi_n+T\psi_k=c_ne^{-a_n^2/2}\eps^n_n(t-a_n)^ne^{(a_n+i\omega_n)t-t^2/2}
+c_ke^{-a_k^2/2}\eps_k^k(t-a_k)^ke^{(a_k+i\omega_k)t-t^2/2}.
$$
On the other hand, $\psi_n+\psi_k$ has only trivial partners, thus there exists
 constants
$c_{k,n}\in\T$, $\eps_{k,n}=\pm1$, $\omega_{k,n},a_{k,n}\in\R$ such that
$$
T(\psi_n+\psi_k)=c_{k,n}e^{-a_{k,n}^2/2}[\eps_{k,n}^nH_n(t-a)+\eps_{k,n}^kH_k(t-a)]
e^{(a_{k,n}+i\omega_{k,n})t-t^2/2}.
$$
Comparing the growth at $\pm\infty$ and $\pm i\infty$ in these two expressions,
we get that the exponential parts have to be the same, that is
$$
a_n+i\omega_n=a_k+i\omega_k=a_{k,n}+i\omega_{k,n}
$$
so that $a_{k,n}=a_n=a_k$ and $\omega_{k,n}=\omega_n=\omega_k$ {\it i.e.} for
every $n$, $a_n=a_0$ and $\omega_n=\omega_0$ as desired. We are then left with
$$
c_n\eps^n_nH_n(t-a_0)+c_k\eps^k_kH_k(t-a_0)=
c_{k,n}\eps_{k,n}^nH_n(t-a_0)+c_{k,n}\eps_{k,n}^kH_k(t-a_0).
$$
But, looking at the highest order term, this implies first that $c_n\eps^n_n=c_{k,n}\eps_{k,n}^n$ and
then $c_k\eps^k_k=c_{k,n}\eps_{k,n}^k$.
If $n$ and $k$ are both even then this reduces further to $c_n=c_{k,n}=c_k$
{\it i.e.} for every $n$ even, $c_n=c_0$.
If $n$ and $k$ are both odd, we get $c_n\eps_n=c_{k,n}\eps_{k,n}=c_k\eps_k$
{\it i.e.} for every $n$ odd, $c_n\eps_n=c_1\eps_1$.
Finally, if $n=0,k=1$ we get $c_0=c_{1,0}$ and
 $c_1\eps_1=c_{1,0}\eps_{1,0}$.
There are thus two alternatives, either $\eps_{1,0}=1$ or $\eps_{1,0}=-1$.
In the first case, $c_1\eps_1=c_0$ and then
$T\psi_n(t)=c_0e^{i\omega_0t}p_n(t-a_0)$. In the second case $c_1\eps_1=-c_0$
so that $c_n\eps_n^n=(-1)^nc_0$ and $T\psi_n(t)=c_0e^{i\omega_0t}\psi_n(-t+a_0)$
as desired.
%
\end{proof}

\section{Pulse type signals}
\label{sec:last}

\subsection{The stability of pulse type signals for the ambiguity
problem} {\,}
\medskip

\noindent The main result in this section can be stated as
follows:

\begin{theorem}
\label{B.pulse}
  Let $0<\eta\leq\frac{1}{3}$ and
$u(t)=\sum_{j=0}^Na_j\chi_{[j,j+\eta]}(t)$ for some
$(a_0,a_1,\ldots,a_N)\in\C^{N+1}$. Then (modulo a trivial
transformation) every solution $v(t)\in L^2(\mathbb{R})$ of the
ambiguity problem (\ref{rap}) is necessarily of the form
$v=\sum_{j=0}^Nb_j\chi_{[j,j+\eta]}$, for some
$(b_0,b_1,\ldots,b_N)\in\C^{N+1}$.
\end{theorem}

This theorem may be seen as an ``uncertainty principle'' for pulse
type signals, in analogy to Lemma \ref{uncert} for Hermite
signals. The techniques we use here, however, are different,
containing ideas from phase retrieval and various limiting
arguments. The role of $\eta\leq\frac13$ is crucial in the proof,
and one may conjecture that $\frac13$ is critical to obtain such
an uncertainty principle.

The following elementary lemma will be used in the sequel.

\begin{lemma}\label{lb0}
Let $u,v$ be Lebesgue measurable functions and $[a,b]\subset\R$.
Assume that for all $x\in[a,b]$, and almost every $t\in\R$,
$u(t)\,v(t+x)=0$. Then, if $t_0\in\supp{u}$ we have $v(t)=0$ for
almost every $t\in t_0+[a,b]$.
\end{lemma}

\begin{proof}
Consider the set
$$
A=\{(t,x)\in\R\times[a,b]\mid u(t)\,v(t+x)\not=0\}.
$$
By Tonelli's theorem and the assumption in the lemma
$$
|A|=\int_{[a,b]}\,\left|\{t\in\R\mid
u(t)\,v(t+x)\not=0\}\right|\,dx=0.
$$
Without loss of generality we shall assume $t_0=0$. For
$0<\eps<\frac{b-a}{2}$, let $U_\eps=\{t\in(-\eps,\eps)\mid
u(t)\not=0\}$ and $V_\eps=\{x\in[a+\eps,b-\eps]\mid v(x)\not=0\}$.
As $0\in\supp{u}$, for every $\eps>0$, $|U_\eps|>0$.

Consider the set $\dst A_\eps = \bigcup_{t\in
U_\eps}\{t\}\times(V_\eps-t)$ and note that
$$
A_\eps \subset \{(t,x)\in U_\eps\times[a,b]\mid
u(t)\,v(t+x)\not=0\}\subset A.
$$
Since $|A|=0$, it follows that $A_\eps$ is measurable in $\R^2$
and $|A_\eps|=0$. Thus, using again Tonelli's theorem
$$
|A_\eps|=\int_{U_\eps}|V_\eps-t|\,dt=|U_\eps|\,|V_\eps|=0.
$$
As $|U_\eps|>0$ this implies that $|V_\eps|=0$ for every $\eps>0$,
thus $|V_0|=0$.
\end{proof}

\begin{proof}[Proof of Theorem \ref{B.pulse}]

\noindent We shall assume $a_0a_N\not=0$. Let $v\in L^2(\R)$ be an
ambiguity partner of $u$, that is
\begin{align}
\left|\mathcal{F}^{-1}\bigl(v\,\overline{v(\cdot-x)}\bigr)(y)\right|=
&
\left|\mathcal{F}^{-1}\bigl(u\,\overline{u(\cdot-x)}\bigr)(y)\right|\notag\\
=& \sum_{k=-N}^N|\mathcal{A} a(k,y)|\,\left|
\frac{\sin(\eta-|x-k|)y/2}{y/2}\right|\,\chi_{[-\eta,\eta]}(x-k),
\label{b.1}
\end{align}
for all $x,y\in\R$. We need to show  that $v$ is a pulse function
of the same type as $u$. This will be obtained directly from
(\ref{b.1}) in various steps. To begin with we recall that, modulo
a trivial transformation, we must have
\begin{equation}
\rm{conv\,}(\supp{v}) = \rm{conv\,}(\supp{u})=[0,N+\eta] \label{b.2}
\end{equation}
(see, e.g., Lemma 1 in \cite[3.2.2]{Ja}). In particular, $v$ is
compactly supported and (\ref{b.1}) is an equality of continuous
functions in $x$ and $y$.

\medskip

{\bf Step 1.} {\sl A bound for the support of $v$.}

\medskip

From (\ref{b.1}) it is clear that, for every
$x\in[\eta,1-\eta]+\Z$,
\begin{equation}
\label{b.3} v(\cdot)\,\overline{v(\cdot+x)}
=u(\cdot)\,\overline{u(\cdot+x)} =0 \ a.e.
\end{equation}
Since $0\in\supp{v}$, we conclude from Lemma \ref{lb0} that
$v(x)=0$, for almost every $x\in[\eta,1-\eta]+\Z$. Thus, there are
some smallest intervals $I_j=[l_j,r_j]\subset j+[-\eta,\eta]$,
$j=0,...,N$, so that
\begin{equation}
\supp{v}\subset\cup_{j=0}^N I_j=\cup_{j=0}^N [l_j,r_j].
\label{b.4}
\end{equation}
Observe that $l_0=0$ and $r_N=N+\eta$ by (\ref{b.2}). Further, we
claim that our assumption $\eta\leq\frac{1}{3}$ actually implies
$r_j-l_j\leq\eta$. Indeed, we already know this for $I_0=
[0,r_0]\subset[0,\eta]$. Let us now show it for $I_1=
[l_1,r_1]\subset[1-\eta,1+\eta]$. Since $l_1\in\supp{v}$, we can
use again Lemma \ref{lb0} and (\ref{b.3}) to conclude
$$
v(l_1+x)=0,\quad\mbox{for almost every }x\in[\eta,1-\eta],
$$
or equivalently, $v$ vanishes in $l_1+[\eta,1-\eta]$. Now, this
interval cannot be strictly contained in $[l_1,r_1]$ because the
latter has length not exceeding $2\eta$ and the former (with left
extreme $l_1+\eta$) has length $1-2\eta\geq\eta$. Therefore, by
the minimality of $I_1$ we must necessarily have $l_1+\eta\geq
r_1$, which gives our claim.  One proceeds similarly with the
other intervals $I_j$.

In particular, we have shown that
$$
\supp{v}\subset\bigcup_{j=0}^N I_j\subset\bigcup_{j=0}^N
[l_j,l_j+\eta].
$$
Observe that we cannot exclude the possibility that some $I_j$ may
be empty. In this case, there is no loss in considering
$l_j=r_j=j$.

\medskip

{\bf Step 2.} {\sl The phase retrieval problem.}

\medskip

Let us now fix $k\in\{0,\ldots,N\}$ and $x\in k+(-\eta,\eta)$. We
then study (\ref{b.2}) as the phase retrieval problem
\begin{eqnarray}
\abs{\ff[v(\cdot)\,\overline{v(\cdot-x)}](y)}&=&
\abs{\ff[u(\cdot)\,\overline{u(\cdot-x)}](y)}\notag\\
& = & \abs{\aa a(k,y)}\,
\left|\frac{\sin(\eta-|x-k|)y/2}{y/2}\right|. \label{b.5}
\end{eqnarray}
By Walther's theorem (\cite{Wa} or \cite[Theorem 2]{Ja}), the
solution to this problem is necessarily of the form
\begin{equation}
\label{b.6} \mathcal{F}[v(\cdot)\,\overline{v(\cdot-x)}](y)=
e^{i\al(x)}\,e^{i\beta(x)y}\, \aa a(k,y)\,
\frac{\sin(\eta-|x-k|)y/2}{y/2}\,G_x(y)
\end{equation}
where $\al(x),\beta(x)$ are real functions, and $G_x$ is a
unimodular function of the form
$$
G_x(y)=\prod_{z\in J_x}\frac{(1-\frac{y}{\overline{z}})\,
e^{\frac{y}{\overline{z}}}}{(1-\frac{y}{z})\,e^{\frac{y}{z}}},
$$
for some set of (non-real) complex numbers $J_x$. The set $J_x$ is
a subset of the {\it complex} zeros of $z\mapsto
\ff[u(\cdot)\,\overline{u(\cdot-x)}](z)$. The effect of $G_x$ is
to take these zeros into their complex conjugates (the so called
{\it zero-flipping}).

Since $z\mapsto\frac{\sin(\eta-|x-k|)z/2}{z/2}$ has only real
zeros, flipping may only occur in the set $\mathcal{Z}_k$ of
\emph{non-real zeros} of $z\mapsto \aa a(k,z)$ (where as usual,
zeros are repeated according to multiplicity). We can partition
$\mathcal{Z}_k=I_x\cup J_x$, with $J_x$ the subset of zeros that
``flip'' in (\ref{b.5}).

Our first claim is that, for each $k=0,\ldots,N$, $J_x$ (and thus
$G_x$) are actually independent of $x\in k+(-\eta,\eta)$. Indeed,
given one such $x_0$ one notices that the holomorphic function
$F_{x_0}(z)= \mathcal{F}[v(\cdot)\,\overline{v(\cdot-x_0)}](z)$ is
not identically zero, and $F_x\to F_{x_0}$ in $\mathcal{H}(\C)$
when $x\to x_0$. Moreover, given any zero
$z\in\mathcal{Z}(F_{x_0})$, by Rouch\'e's theorem we obtain the
equality of multiplicities $m(z,F_x)=m(z,F_{x_0})$, for all
$|x-x_0|<\eps$ provided $\eps=\eps(x_0,z)>0$ is small enough.

Proceeding as before for every $x_0$, an easy
compactness-connectedness argument gives $m(z,F_x)=m(z,F_k)$ for
all $x\in k+(-\eta,\eta)$. Finally, repeating this argument with
all zeros $z\in\mathcal{Z}(F_k)$ one concludes $J_x=J_k$ for all
$x\in k+(-\eta,\eta)$. We will then write $G_k=G_x$ for such $x$.

\medskip

{\bf Step 3.} {\sl Determination of the support of
$v(\cdot)\,\overline{(\cdot-x)}$.}

\medskip

Let us now go back to (\ref{b.6}) and define the bounded function
\begin{equation}
{\widehat U}_x(y)= e^{i\al(x)}\,e^{i\beta(x)y}\, \aa a(k,y)\,
G_k(y), \label{b.7}
\end{equation}
so that $U_x$ is a tempered distribution satisfying, for all $x\in
k+(-\eta,\eta)$,
\begin{equation}
v(\cdot)\,\overline{v(\cdot-x)}=
U_x*\chi_{[-\frac{\eta-|x-k|}2,\frac{\eta-|x-k|}2]}. \label{b.8}
\end{equation}
Next, we define another distribution ${\widetilde U}_k$ by
$$
\widehat{{\widetilde U}}_k(y)=\aa a(k,y)\, G_k(y),
$$
so that, for $x\in k+(-\eta,\eta)$,
\begin{equation}
{\widetilde U}_k=e^{-i\al(x)}\,U_x(\cdot+\beta(x)). \label{b.9}
\end{equation}
Let us emphasize that, in this identity, ${\widetilde U}_k$ does
not depend on $x$. Now, if we consider $k=0$ and fix
$x\in[0,\eta)$ we must have, using step 1,
\begin{eqnarray}
\cup_{j=0}^N [l_j+x,l_j+\eta]& \supset &
\supp{v(\cdot)\,\overline{v(\cdot-x)}}\nonumber\\
 & = & \supp{U_x*\chi_{[-\frac{\eta-x}2,\frac{\eta-x}2]}},\nonumber\\
& = & \supp{\left({\widetilde
U}_0*\chi_{[-\frac{\eta-x}2,\frac{\eta-x}2]}\right)}+\beta(x).
\label{b.10}
\end{eqnarray}
Now, as $v(\cdot)\,\overline{v(\cdot-x)}$ is supported in
$[0,N+\eta]$, $z\mapsto
\mathcal{F}[v(\cdot)\,\overline{v(\cdot-x)}](z)$ is entire of
exponential type at most $N+\eta$ (for any $x$). It follows that
$\beta$ is a bounded function and thus we may find a sequence
$x_m\nearrow\eta$ so that $\beta(x_m)$ has some limit, say
$\beta_+$.

Next recall the following elementary fact:

\noindent{\sl for every distribution $U\in\mathcal{S}'$ we have
$$
\frac{1}{2\delta}\,U*\chi_{-(\delta,\delta)}\to U
$$
when $\delta\to0$ with convergence in $\mathcal{S}'$.}

Then, letting $x_m\to\eta$ in (\ref{b.10}) we easily obtain
$$
\supp{{\widetilde U}_0}\subset \cup_{j=0}^N\{l_j+\eta\} -\beta_+.
$$
Further, observe that $\widehat{\widetilde{U}}_0$ is bounded (and
hence cannot be a polynomial), which necessarily implies
\begin{equation}
{\widetilde U}_0=\sum_{j=0}^N\gamma_j\dt_{l_j+\eta-\beta_+},
\label{b.11}
\end{equation}
for some complex numbers $\gamma_j$, $j=0,\ldots,N$. Thus, we
conclude that, if $x\in(-\eta,\eta)$,
\begin{equation}
U_x= e^{i\al(x)}
\sum_{j=0}^N\gamma_j\dt_{l_j+\eta+\beta(x)-\beta},\label{b.12}
\end{equation}
and therefore
\begin{equation}
v(\cdot)\,\overline{v(\cdot-x)}=
 e^{i\al(x)}\sum_{j=0}^N\gamma_j
\chi_{[-\frac{\eta-|x|}2,\frac{\eta-|x|}2]+l_j+\eta+\beta(x)-\beta_+}.
\label{b.13}
\end{equation}

\medskip

{\bf Step 4.} {\sl Determination of $|v|$.}

\medskip

We begin by showing that $\gamma_0\gamma_N\not=0$. Indeed, we test
(\ref{b.13}) with $x=0$, and using the property that
$0\in\supp{v}$, we find a smallest integer $j_0\in\{0,\ldots,N\}$
such that
$$
0\in [-\frac{\eta}{2},\frac{\eta}{2}]+l_{j_0}+\eta+\beta(0)-\beta.
$$
We claim that $j_0=0$. If not we must have
$\frac{3\eta}{2}+\beta(0)-\beta<0$ (since $l_0=0$), and thus
$\beta(0)-\beta<-\frac{3\eta}{2}$. But now, since
$N+\eta\in\supp{v}$ we also have
$$
N+\eta\leq\frac{\eta}{2}+l_N+\eta+\beta(0)-\beta<
\tfrac{3\eta}{2}+l_N-\frac{3\eta}{2}=l_N,
$$
which is a contradiction (since $r_N=N+\eta\in\supp{v}$). Thus,
$j_0=0$, which forces $\gamma_0\not=0$. A completely symmetrical
argument gives $\gamma_N\not=0$.

Next, we shall determine explicitly the function $\beta(x)$ in
(\ref{b.6}). Recall from (\ref{b.10})  that
$$
\supp{v(\cdot)\,\overline{v(\cdot-x)}}\subset
\begin{cases}
\bigcup_{j=0}^N[x+l_j,l_j+\eta], &\mbox{if }x\in[0,\eta)\\
\bigcup_{j=0}^N[l_j,l_j+\eta+x], &\mbox{if }x\in(-\eta,0]\\
\end{cases}.
$$
Since $\gamma_0\gamma_N\not=0$, we see from the (\ref{b.13}) that
the extreme points
$-\frac{\eta-|x|}{2}+l_0+\eta+\beta(x)-\beta_+$,
$\frac{\eta-|x|}{2}+l_N+\eta+\beta(x)-\beta_+$ must belong to
$\supp{v(\cdot)\,\overline{v(\cdot-x)}}$. Therefore, if
$x\in[0,\eta)$,
$x+l_0\leq-\frac{\eta-x}{2}+l_0+\eta+\beta(x)-\beta_+$ so that
$$
-\frac{\eta-x}{2}\leq \beta(x)-\beta_+,
$$
and $\frac{\eta-|x|}{2}+l_N+\eta+\beta(x)-\beta_+\leq l_N+\eta$ so
that
$$
\beta(x)-\beta_+\leq -\frac{\eta-x}{2}.
$$
Thus, we conclude $\beta(x)=\beta_+-\frac{\eta-x}{2}$,
$x\in[0,\eta)$. Proceeding symmetrically with $x\in(-\eta,0]$ one
extends this identity to all $x\in(-\eta,\eta)$. In conclusion,
going back to equation (\ref{b.13}) with $x=0$ we have shown that
$$
 |v|^2= e^{i\al(0)}
\sum_{j=0}^N\gamma_j\,\chi_{[l_j,l_j+\eta]}.
$$
Next we shall determine  explicitly the values of $l_j$. As we
said in step 1, there is no loss in assuming $l_j=j$ when
$\gamma_j=0$. We will prove that we must also have  $l_j=j$ when
$\gamma_j\not=0$. Indeed, we already know that $l_0=0$. Moreover,
when $\gamma_j\not=0$ we know from step 1 that
$$
[l_j,l_j+\eta]\,\subset\,j+[-\eta,\eta],
$$
from which it follows $j-\eta\leq l_j\leq j$. Assume by
contradiction that for one such $j$ we have $j-\eta\leq l_j\leq
j-\eps$, for some $0<\eps<\eta$. Then we can select
$$
x=l_j-(\eta-\eps)\in(j-1)+[\eta,1-\eta],
$$
so that by (\ref{b.3}) it holds $v(\cdot)\overline{v(\cdot-x)}=0$,
$a.e.$ Now, when $t\in[l_j,l_j+\eps]$ we also have
$t-x\in[\eta-\eps,\eta]\subset[0,\eta]$, and therefore, for
$t\in[l_j,l_j+\eps]$,
$$
v(t)\,\overline{v(t-x)}=\gamma_j\,{\overline{\gamma_0}}\not=0,
$$
which is a contradiction. Thus, we have proven
\begin{equation}
|v|^2= e^{i\alpha(0)} \sum_{j=0}^N\gamma_j\chi_{[j,j+\eta]},
\label{b.15}
\end{equation}
or more generally, looking at (\ref{b.13}), for $x\in(-\eta,\eta)$
\begin{equation}
v(\cdot)\,\overline{v(\cdot-x)}=
 e^{i\alpha(x)}
\sum_{j=0}^N\gamma_j\chi_{\ent{-\frac{\eta-|x|}2,\frac{\eta-|x|}2}+j+\frac{x+\eta}2}.
 \label{b.14}
\end{equation}

\medskip

{\bf Step 5.} {\sl Determination of the phase of $v$.}

\medskip

From (\ref{b.15}) we conclude that there are numbers
$b_0,\ldots,b_N\geq0$, and a function $t\mapsto\phi(t)$ real such
that
\[
v(t)=
 e^{i\phi(t)}
\sum_{j=0}^Nb_j\,\chi_{[j,j+\eta]}.
\]
Observe that we can modify $v$ in a null set so that this equality
holds in all points $t\in\R$. We want to show that the phase
$\phi(t)$ is constant in each interval $[j,j+\eta]$ for which
$b_j\not=0$. When $x\in[0,\eta)$, using the expression in
(\ref{b.14}) we see that
$$
v(\cdot)\,\overline{v(\cdot-x)}=
e^{i\al(x)}\sum_{j=0}^N\gamma_j\,\chi_{[x,\eta]+j}
=e^{i(\phi(t)-\phi(t-x))}\sum_{j=0}^N|b_j|^2\,\chi_{[x,\eta]+j}.
$$
Since by (\ref{b.15}) $e^{i\al(0)}\gamma_j\geq0$, we must have
$$
\widetilde{\alpha}(x)\equiv\al(x)-\al(0)=\phi(t)-\phi(t-x),\quad(\mbox{mod}{2\pi})
$$
whenever $x\in[0,\eta)$, $t\in[x,\eta]+j$ and $b_j\not=0$.
Choosing $t=x+j$ we see that
$\phi(x+j)=\widetilde{\alpha}(x)+\phi(j)$ (mod $2\pi$), and
therefore
\[
\widetilde{\alpha}(x)=\widetilde{\alpha}(t)-\widetilde{\alpha}(t-x)\quad(\mbox{mod}{2\pi}),\quad
x,t,\,t-x\in[0,\eta).
\]
This is equivalent to
\[
\widetilde{\alpha}(t+x)=\widetilde{\alpha}(t)+\widetilde{\alpha}(x)\quad(\mbox{mod}{2\pi}),\quad\mbox{when}\quad
x,t,t+x\in[0,\eta),
\]
which by continuity of $\widetilde{\alpha}$ (by (\ref{b.6}))
implies $\widetilde{\alpha}(t)=\omega t$, $t\in[0,\eta)$, for some
real number $\omega$. Thus, modulo $2\pi$,
$\phi(t+j)=\phi(j)+\omega t$, $t\in[0,\eta)$, so calling
$\widetilde{b}_j=e^{i(\phi(j)-\omega j)}b_j$ we conclude
\[
v(t)=e^{i\omega
t}\,\sum_{j=0}^N\widetilde{b}_j\,\chi_{[j,j+\eta]}.
\]
Therefore we have shown that, modulo a trivial transformation, $v$
is a signal of pulse type of the same form as $u$, concluding the
proof of the theorem.
\end{proof}

\subsection{Rareness of pulse signals with non-trivial partners}

Contrary to section \ref{sec:3}, from now on it will be more convenient to
study the Discrete Radar Ambiguity Problem (P) for sequences rather than
the Periodic Radar Ambiguity Problem ($\widehat{P}$). Let us first note that there is no difficulty
to transpose Proposition \ref{tl2} to this context. Note that the trivial solutions are
generated by the two representations of the periodized Heisenberg
group $\hil=\T\times\T\times\Z$ on $\ell^2(\Z)$ given as follows.
For $h=(\beta,\omega,l)\in\hil$ and
$a=(a_j)_{j\in\Z}\in\ell^2(\Z)$, define $b=S_ha$ by
$$
b_j=e^{i\beta+ij\omega}a_{j-l},
$$
and $\widetilde{b}=\widetilde{S}_ha$ by
$$
\widetilde{b}_j=e^{i\beta+ij\omega}a_{-j-l}.
$$
Further, when looking for partners of a finite sequence $a$, we may replace $a$ by a trivial partner
and assume that $a=(a_0,\ldots,a_N)$ for some integer $N$ and that $a_0a_N\not=0$. We will then
write $a\in\ss(N)$. Transposing Lemma \ref{tl3} from trigonometric polynomials to finite
sequences, a partner of $b$ of $a$ may then also be assumed to be in $\ss(N)$.

In view of Theorem \ref{B.pulse}, the study of problem \ref{rap} for pulse type signals of finite length
is then reduced to the following finite dimensional ambiguity problem, where $N$ is a fixed positive integer.

\medskip

\noindent {\bf Ambiguity problem in $\ss(N)$.}
{\sl Given  $a=(a_0,a_1,\ldots, a_N)\in\ss(N)$, find all
$b=(b_0,b_1,\ldots, b_N)\in\ss(N)$ such that
\begin{equation}
\label{eqdiscret}
\abs{\aa(b)(j,y)}=\abs{\aa(a)(j,y)}\qquad\mbox{for all }j\in\Z,\ y\in\T.
\end{equation}}

\medskip

We will now use the following notation.
\begin{notation}
If $b$ is a trivial ambiguity partner of $a$, we write $b\tdrap a$.
If $b\drap a$ but $b\not\tdrap a$, we call $b$ a {\it strange partner} of $a$
and write $b\ntdrap a$.
\end{notation}

The goal is to describe the class of all signals $a$ which only
admit trivial partners $b$. Several results in this direction have already appeared in
\cite{GJP}, which we describe now. We shall denote the complementary of the searched class by
$$
\ee(N)=\set{a\in\ss(N)\,:\ a \mbox{ admits strange partners}}.
$$
It is easy to see that $\ee(N)=\emptyset$ for $N=0,1,2$. The main
result in \cite{GJP} establishes that for larger values of $N$
this set cannot be too large.

\begin{theorem}
\label{th:pulserare}
For every $N\geq 3$, $\ee(N)$ is a non-empty semi-algebraic
variety of real dimension at most $2N+1$.
\end{theorem}

We recall that a semi-algebraic variety is a set defined by
polynomial equalities and/or inequalities. The theorem says that
$\ee(N)$ has this structure, and moreover is contained in a real
algebraic variety (i.e., finite unions of polynomial zero sets) of real
dimension $2N+1$. This implies that $\ee(N)$ has Lebesgue measure
$0$ in $\C^{N+1}$ and is also thin in the Baire sense.

\begin{corollary}
For every $N\geq 0$, quasi-all and almost all elements of $\ss(N)$ have only
trivial partners.
\end{corollary}

A full description of $\ee(N)$ for  $N=3,4$ can be found in
\cite{GJP}. In particular, $\ee(3)$ contains sequences with all
$a_j\not=0$, $j=0,1,2,3$. This shows that sequences with strange
partners do not necessarily have to contain ``gaps'', a remarkable
fact in view of the results in Section \ref{sec:3}. In \cite{GJP}, a general
argument showing the non-emptiness of $\ee(N)$ for $N\geq3$ was
only sketched. The object of the next section is to prove it in
full detail.


\subsection{Construction of strange partners}
\label{sec:4.3}

A simple way to construct strange ambiguity partners when $N=2K+1$
is odd is as follows: take $\alpha=(\alpha_0,\ldots,\alpha_K)$ be
any sequence of length $K$. A direct computation of their
ambiguity functions shows that for $\lambda\in\C$, the sequences
\begin{equation}
a_k=\begin{cases}
\alpha_p&\mbox{when }k=2p\\
\lambda\alpha_p&\mbox{when }k=2p+1\\
\end{cases}
\quad\mbox{and}\quad b_k=\begin{cases}
\overline{\lambda}\alpha_p&\mbox{when }k=2p\\
\alpha_p&\mbox{when }k=2p+1\\
\end{cases}
\label{ex:triv}
\end{equation}
are ambiguity partners. In general, these are non-trivial partners
(see \cite[p. 102]{GJP}). Since this method is restricted to $N$
odd, we will now describe another method that gives elements of
$\ee(N)$ as soon as $N\geq 4$.

First recall from \cite{GJP} that when $a\in\ss(N)$ one can
reformulate (\ref{eqdiscret}) as an equivalent combinatorial
problem on matrices. Namely, if we let $K_a$ be the matrix with
entries
$$
d_{j,k}=\begin{cases} a_{\frac{j+k}{2}}a_{\frac{j-k}{2}}&\mbox{if
}j,k\mbox{
have same parity}\\
0&\mbox{else}\\ \end{cases},
$$
then we have the following

\begin{proposition} Two sequences $a,b\in\ss(N)$ are ambiguity partners if and
only if
$$
K_a^*K_a=K_b^*K_b.
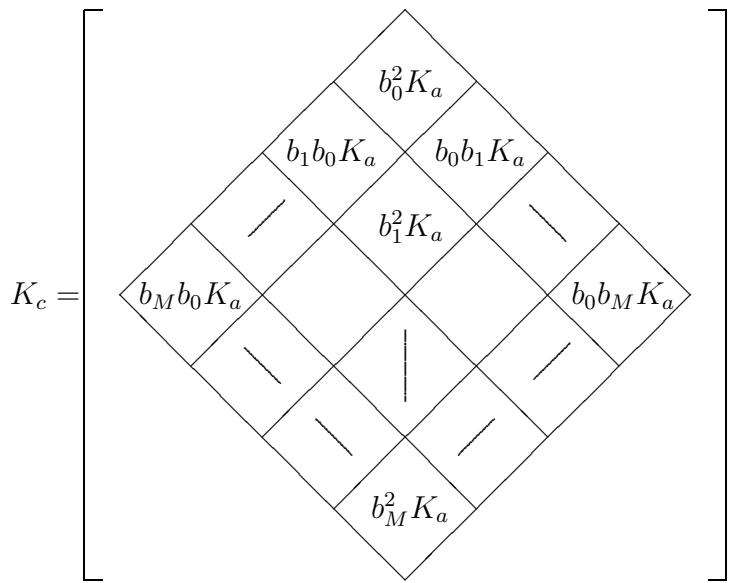
$$
\end{proposition}

\begin{example}
If $a\in\ss(5)$, the matrix of $K_a$ is given by
$$
\ka{a_0}{a_1}{a_2}{a_3}{a_4}{a_5}
$$
(non written elements of that matrix are $0$).
\end{example}


We shall make use of the \emph{Kronecker product} of matrices,
which for $A$ and $B=[b_{i,j}]_{1\leq i,j\leq n}$ is the matrix
defined by blocks as
$$
A\otimes B=\ent{\begin{matrix}
Ab_{1,1}&Ab_{1,2}&\ldots&Ab_{1,n}\\
Ab_{2,1}&Ab_{2,2}&\ldots&Ab_{2,n}\\
\vdots&\vdots&\ddots&\vdots\\
Ab_{n,1}&Ab_{n,2}&\ldots&Ab_{n,n}\\
\end{matrix}}.
$$
This product has the following elementary properties:
\begin{itemize}
 \item[---] $(A\otimes B)^*=A^*\otimes B^*$, \item[---]
$(A\otimes B)(C\otimes D)=(AC)\otimes(BD)$.
\end{itemize}

We shall compute the Kronecker product of two ambiguity matrices
$K_a$ and $K_b$ and show that it corresponds to the ambiguity
matrix of a new sequence $c$ produced by a certain product rule
involving $a$ and $b$. This turns out to produce many natural
examples of sequences with strange partners.

For this, it is convenient to change the way to enumerate the
entries of such matrices, by introducing the following
``lattice coordinates'': let $\gamma=\begin{bmatrix}-1&1\\ 1&1\\
\end{bmatrix}$ and $\Gamma=\gamma\Z^2$ be a sub-lattice of $\Z^2$.
Given $N\geq1$ we consider the subset of entries
$\Gamma_N=\left\{\begin{bmatrix}-1&1\\ 1&1\\
\end{bmatrix}\begin{bmatrix}m\\ l\\
\end{bmatrix}\ :\ 0\leq m,l\leq N\right\}$. If
$a=(a_0,a_1,\ldots, a_N)\in\mathbb{C}^{N+1}$, then $K_a$ is
supported in $\Gamma_N$ and
$$
\left(K_a\right)_{i,j}=a_ma_\ell\,\mbox{ if }
\begin{bmatrix}i\\ j\\ \end{bmatrix}=\begin{bmatrix}-1&1\\ 1&1\\
\end{bmatrix}
\begin{bmatrix}m\\\ell\\ \end{bmatrix}\ :\ 0\leq m,\ell\leq N.
$$
Thus, $K_a$ is completely determined by the matrix
$\tK_a[m,\ell]:=\left(K_a\right)_{i,j}$ when $\begin{bmatrix}i\\ j\\
\end{bmatrix}=\gamma
\begin{bmatrix}m\\ \ell\\ \end{bmatrix}$.

\begin{lemma} \label{kaokb}
Let $a=(a_0,\ldots,a_N)$ and $b=(b_0,\ldots,b_M)$ be two finite
sequences with associated polynomials
$P(z)=\dst\sum_{k=0}^Na_kz^k$, $Q(z)=\dst\sum_{k=0}^Mb_kz^k$.
Consider the polynomial
$P(z)Q(z^{N+1})=\dst\sum_{k=0}^Kc_kz^k$ 
and let $c=(c_0,\ldots,c_K)$. Then the ambiguity matrix $K_c$ is
supported in $\Gamma_N+(N+1)\Gamma_M$ and satisfies
\begin{equation}\label{tKa}
    \tK_c[i+(N+1)m,j+(N+1)\ell] =b_mb_\ell a_ia_j.
\end{equation}
In particular, $\tK_c=\tK_a\otimes\tK_b$ and the matrix $K_c$ can
be drawn as
\begin{center}
\setlength{\unitlength}{2947sp}%
\begingroup\makeatletter\ifx\SetFigFont\undefined%
\gdef\SetFigFont#1#2#3#4#5{%
  \reset@font\fontsize{#1}{#2pt}%
  \fontfamily{#3}\fontseries{#4}\fontshape{#5}%
  \selectfont}%
\fi\endgroup%
\begin{picture}(5937,4824)(676,-4273)
\thinlines
{\put(1501,-1861){\line( 1, 1){2400}}
\put(3901,539){\line( 1,-1){2400}}
\put(6301,-1861){\line(-1,-1){2400}}
\put(3901,-4261){\line(-1, 1){2400}}
}%
{\put(2101,-1261){\line( 1,-1){2400}}
}%
{\put(2701,-661){\line( 1,-1){2400}}
}%
{\put(3301,-61){\line( 1,-1){2400}}
}%
{\put(4501,-61){\line(-1,-1){2400}}
}%
{\put(5101,-661){\line(-1,-1){2400}}
}%
{\put(5701,-1261){\line(-1,-1){2400}}
}%
{\put(1351,539){\line(-1, 0){150}}
\put(1201,539){\line( 0,-1){4800}}
\put(1201,-4261){\line( 1, 0){150}}
}%
{\put(6451,-4261){\line( 1, 0){150}}
\put(6601,-4261){\line( 0, 1){4800}}
\put(6601,539){\line(-1, 0){150}}
}%
{\multiput(3901,-2161)(0.00000,-8.95522){68}{\makebox(1.6667,11.6667){\SetFigFont{5}{6}{\rmdefault}{\mddefault}{\updefault}.}}
}%
{\multiput(2551,-2311)(6.38298,-6.38298){48}{\makebox(1.6667,11.6667){\SetFigFont{5}{6}{\rmdefault}{\mddefault}{\updefault}.}}
}%
{\multiput(3151,-2911)(6.38298,-6.38298){48}{\makebox(1.6667,11.6667){\SetFigFont{5}{6}{\rmdefault}{\mddefault}{\updefault}.}}
}%
{\multiput(4951,-1111)(6.38298,-6.38298){48}{\makebox(1.6667,11.6667){\SetFigFont{5}{6}{\rmdefault}{\mddefault}{\updefault}.}}
}%
{\multiput(4651,-2911)(-6.38298,-6.38298){48}{\makebox(1.6667,11.6667){\SetFigFont{5}{6}{\rmdefault}{\mddefault}{\updefault}.}}
}%
{\multiput(5288,-2273)(-6.38298,-6.38298){48}{\makebox(1.6667,11.6667){\SetFigFont{5}{6}{\rmdefault}{\mddefault}{\updefault}.}}
}%
{\multiput(2888,-1073)(-6.38298,-6.38298){48}{\makebox(1.6667,11.6667){\SetFigFont{5}{6}{\rmdefault}{\mddefault}{\updefault}.}}
}%
\put(576,-1936){\makebox(0,0)[lb]{\smash{\SetFigFont{12}{14.4}{\rmdefault}{\mddefault}{\updefault}{$K_c=$}%
}}}
\put(3676,-151){\makebox(0,0)[lb]{\smash{\SetFigFont{12}{14.4}{\rmdefault}{\mddefault}{\updefault}{$b_0^2K_a$}%
}}}
\put(3616,-3736){\makebox(0,0)[lb]{\smash{\SetFigFont{12}{14.4}{\rmdefault}{\mddefault}{\updefault}{$b_M^2K_a$}%
}}}
\put(3676,-1361){\makebox(0,0)[lb]{\smash{\SetFigFont{12}{14.4}{\rmdefault}{\mddefault}{\updefault}{$b_1^2K_a$}%
}}}
\put(1661,-1936){\makebox(0,0)[lb]{\smash{\SetFigFont{12}{14.4}{\rmdefault}{\mddefault}{\updefault}{$b_Mb_0K_a$}%
}}}
\put(2901,-736){\makebox(0,0)[lb]{\smash{\SetFigFont{12}{14.4}{\rmdefault}{\mddefault}{\updefault}{$b_1b_0K_a$}%
}}}
\put(4151,-736){\makebox(0,0)[lb]{\smash{\SetFigFont{12}{14.4}{\rmdefault}{\mddefault}{\updefault}{$b_0b_1K_a$}%
}}}
\put(5301,-1936){\makebox(0,0)[lb]{\smash{\SetFigFont{12}{14.4}{\rmdefault}{\mddefault}{\updefault}{$b_0b_MK_a$}%
}}}
\end{picture}
.\end{center} \noindent
\end{lemma}

\begin{proof}[Proof] As noted
before
$$
\tK_c[i+(N+1)m,j+(N+1)\ell]=c_{i+(N+1)m}\,c_{j+(N+1)\ell}.
$$
Now by construction of $c$, the only non-null coefficients are
$c_{i+(N+1)m}=a_ib_m$ for $0\leq i\leq N$ and $0\leq m\leq M$.
This gives (\ref{tKa}). To justify the drawing observe that the
submatrix with coordinates in $\Gamma_N +(N+1)\gamma
\begin{bmatrix}m\\ \ell\\ \end{bmatrix}$ is precisely
$b_mb_\ell K_a|_{\Gamma_N}$, which as $m$ and $\ell$ moves fills
each of the parallelograms in the picture.
\end{proof}

\begin{lemma}
\label{kr2} Let $a=(a_0,\ldots,a_N)$ and $b=(b_0,\ldots,b_M)$ be
two finite sequences with associated polynomials
$P(z)=\dst\sum_{k=0}^Na_kz^k$, $Q(z)=\dst\sum_{k=0}^Mb_kz^k$.
Consider this time the polynomial
$P(z)Q(z^{2N+1})=\dst\sum_{k=0}^Kc_kz^k$  and let
$c=(c_0,\ldots,c_K)$.
Then the ambiguity matrix of $c$ is $K_c=K_a\otimes K_b$.
\end{lemma}

\begin{proof}[Proof] Applying the previous lemma to
$\widetilde{P}(z)=\dst\sum_{k=0}^{2N}\widetilde{a_k}z^k$ where
 $\dst\begin{cases}
\widetilde{a_i}=a_i&\text{if }0\leq i\leq N\\
\widetilde{a_i}=0&\text{if }N\leq i\leq2N\\ \end{cases}$ we see
that \[ \supp K_{c}\subset \supp K_{\tilde a} + (2N+1)\supp K_{
b}.
\] Since $K_{\tilde a}$ vanishes in $\Gamma_{2N}\setminus\Gamma_N$,
we actually have
\[ \supp K_{ c}\subset \supp K_{a} + (2N+1)\supp
K_{ b}.
\]
If we regard $K_a$ as a square (2N+1)-matrix, this implies that
$K_c$ can be written as a collection of disjoint consecutive
square blocks
$\{K_a+(2N+1)\begin{bmatrix}i\\ j\\
\end{bmatrix}\;\colon\;\begin{bmatrix}i\\ j\\
\end{bmatrix}\in\supp K_b\}$.
Next, if we take  $\begin{bmatrix}i\\ j\\
\end{bmatrix}=\gamma
\begin{bmatrix}m\\ \ell\\ \end{bmatrix}\in\supp
K_b\subset\Gamma_M$, then by the previous lemma the value of $K_c$
in the corresponding block is precisely
\[
b_mb_\ell K_a|_{\Gamma_N}.
\]
This shows $K_c=K_a\otimes K_b$ as asserted.
\end{proof}

A sequence $c$ constructed from $a$ and $b$ as in the statement of
Lemma \ref{kr2} will be denoted by $c=a\otimes b$. Recall also
that $a\simeq b$ means that $a$ and $b$ are ambiguity partners as
in (\ref{eqdiscret}).

\begin{cor} Let $a,b,a',b'$ be four finite sequences.
If $a\drap a'$ and $b\drap b'$, then $a\otimes b\drap a'\otimes
b'$.
\end{cor}

\begin{proof}[Proof] From the previous lemma and elementary properties of the Kronecker
product we see that
\begin{align}
K_{a'\otimes b'}^*K_{a'\otimes b'}^{}=&(K_{a'}\otimes K_{b'})^*
(K_{a'}\otimes K_{b'})\notag\\
=&(K_{a'}^*\otimes K_{b'}^*)(K_{a'}^{}\otimes
K_{b'}^{})=(K_{a'}^*K_{a'}^{})\otimes(K_{b'}^*K_{b'}^{}\notag)\\
=&(K_a^*K_a^{})\otimes(K_b^*K_b^{})=K_{a\otimes b}^*K_{a\otimes
b}^{}.\notag
\end{align}
Thus, $a\otimes b\drap a'\otimes b'$ as asserted.
\end{proof}

This corollary enables us to construct sequences $a\in\sN$ with
strange partners, as soon as $N\geq4$.

\begin{example} Let $a=(1,2)$, $b=(1,2)$ and $b'=(2,1)$, then $a\otimes
b\drap a\otimes b'$. But
$$
a\otimes b=(1,2,0,2,4) \quad\mathrm{whereas}\quad a\otimes
b'=(2,4,0,1,2)
$$
so that  $a\otimes b$ and $a\otimes b'$ are not trivial partners
and $a\otimes b\in\ee(4)$. Moreover, applying the above
construction to $(1,2,0,\ldots,0)$ regarded as sequence in
$\C^{N+1}$, we obtain a sequence $a\otimes b\in\ss(2N+2)$, which
shows that $\ee(2N+2)\not=\emptyset$ for all $N \geq 1$.
\end{example}

\begin{example}
Other examples can be produced by iterating this process. For
instance, consider the sequence $c$ associated with the polynomial
\[
R(z)=\prod_{j=0}^J (\alpha_j+ \beta_jz^{3^j}).
\]
Non-trivial ambiguity partners can be obtained by selecting a
collection of $j$'s and replacing the corresponding factors in the
polynomial by $\alpha_j+ c_j\beta_jz^{3^j}$ or $\beta_j+
c_j\alpha_jz^{3^j}$, with $|c_j|=1$. It is possible to show
(although harder) that these are all the possible ambiguity
partners of $c$. Observe finally that these kind of examples are
of a different nature than those in Proposition \ref{prop1}.
\end{example}

As an application we obtain the following remarkable result.

\begin{cor}
The set of all functions $u\in L^2(\R)$ having strange ambiguity
partners in the sense of (\ref{rap}) is dense in $L^2(\R)$.
\end{cor}

\begin{proof}[Proof] Let $f\in L^2(\R)$, which we may assume
with $\norm{f}\leq 1$. Given $0<\eps<1$ we can find $f_c$ with
compact support such that $\norm{f-f_c}<\eps$. Suppose that $\supp
f_c\subset[-\frac R4,\frac R4]$.

Further, taking $a=(1,\eps)$, $b=(1,\eps)$ and $b'=(\eps,1)$, then
$a\otimes b\drap a\otimes b'$. But
$$
a\otimes b=(1,\eps,0,\eps,\eps^2) \quad\mathrm{whereas}\quad
a\otimes b'=(\eps,\eps^2,0,1,\eps)
$$
so that the pulse type signals
$$
u(t)=\sum (a\otimes b)_jf_c(t-Rj) \mbox{ and } v(t)=\sum (c\otimes
b)_jf_c(t-Rj)
$$
are non-trivial ambiguity partners and
$$
\norm{f-u}\leq \norm{f-f_c}+(2\eps+\eps^2)\norm{f_c}\leq 7\eps.
$$
\end{proof}

\section{Conclusion}

The radar ambiguity problem is a difficult and still widely open
problem. In this paper we have concentrated in the most common
classes of signals (Gaussian and rectangular pulses), and shown
how to tackle such cases with real and complex analysis methods,
and also with algebraic approaches. We are still unable to say
much about the general case, but the originality of our methods
may be useful when studying similar problems in the phase
retrieval literature.

For Hermite functions, we rediscover a conjecture from the 70's
which is stronger than the uncertainty principle for ambiguity
functions in Section 2.1. We are almost certain that Hermite
functions must have only trivial partners.
 Indeed, we have only used a
small part of the relations between partners to conclude in the
generic case. On the other side, our proof becomes technically
very complicate when dealing with other cases, and new  ideas may
be necessary.

In the case of pulse type signals, we have both the rareness of
functions with strange partners, some criteria to have only
trivial solutions ({\it see} \cite{GJP}) and various ways to
construct functions that have strange partners. On the other hand,
we are unable to attack the discrete problem, that is Problem
(\ref{pbdisc}), for general sequences with infinite length. We
know that sequences with strange partners are dense (as well as
those with only trivial partners), but it seems likely to us that
they must be ''small" in a suitable sense (such as Baire
category), although we still lack of evidence for this.

Let us conclude by saying that more general classes would be of interest
for instance compactly supported functions ({\it see} \cite{Ja} for some results)
and functions of the form $P(x)e^{-x^2/2}$ with $P$ an entire function of order $<1$.
For the later, note that our techniques do not allow to say anything since
we always start with the highest order coefficient of $P$ when $P$ is a polynomial
(it may be shown that every ambiguity partner is of the same form).

\bibliographystyle{plain}

\end{document}